\documentclass[11pt]{article}
\usepackage{amsfonts,latexsym,amsmath,amscd,geometry}
\usepackage{amssymb}
\usepackage{latexsym}
\usepackage{xcolor}
\usepackage[title]{appendix}

\newcommand \nc{\newcommand}

\numberwithin{equation}{section}
\newtheorem{theorem}{Theorem}[section]
\newtheorem{lemma}[theorem]{Lemma}
\newtheorem{proposition}[theorem]{Proposition}

\newtheorem{definition}[theorem]{Definition}

\newtheorem{remark}[theorem]{Remark}

\nc{\ba}{\begin{array}}\nc{\ea}{\end{array}}
\nc{\be}{\begin{eqnarray}}\nc{\ee}{\end{eqnarray}}
\nc{\beq}{\begin{equation}}\nc{\eeq}{\end{equation}}
\nc{\bex}{\begin{eqnarray*}}\nc{\eex}{\end{eqnarray*}}
\nc{\btm}{\begin{theorem}} \nc{\etm}{\end{theorem}}
\nc{\blm}{\begin{lemma}} \nc{\elm}{\end{lemma}}
\nc{\R}{\mathbb{R}}  
\def\x{\mathbf{x}}\def\e{\mathbf{e}}
\def\T{\mathbf{T}}
\def\pf{\noindent{\bf Proof.\quad}}\def\endpf{\hfill$\Box$}
\def\di{\mbox{div\,}}
\def\curl{\mbox{curl\,}}

\def\u{\mathbf{u}}
\def\d{\mathbf{d}}

\def\v{\varepsilon}

\begin{document}

\title{Poiseuille flow of hyperbolic Ericksen-Leslie system in dimension two}

\author{Geng Chen
\footnote{Department of Mathematics, University of Kansas, Lawrence, KS 66045, U.S.A. Email: gengchen@ku.edu}\quad
Tao Huang
\footnote{Department of Mathematics, Wayne State University, Detroit, MI, 48202, U.S.A. Email: taohuang@wayne.edu}\quad 
Xiang Xu
\footnote{Department of Mathematics, University of Kansas, Lawrence, KS 66045, U.S.A. Email: x354x351@ku.edu}\quad 
Qingtian Zhang
\footnote{School of Mathematical Sciences, Shenzhen University, Shenzhen 518060, Guangdong, P.R.C. Email: zhang$\_$qingtian@szu.edu.cn}
}
\maketitle

\begin{abstract} 
In this paper, we study the Poiseuille laminar flow in a tube for the full Ericksen-Leslie system. It is a parabolic-hyperbolic coupled system which may develop singularity in finite time. We will prove the global existence of energy weak solution, and the partial regularity of solution to  system. We first construct global weak finite energy solutions by the Ginzburg-Landau approximation and the fixed-point arguments. Then we obtain the enhanced regularity of solution. Different from the solution in one space dimension, the finite energy solution of Poiseuille laminar flow in a tube may still form a discontinuity at the origin. We show that at the first possible blowup time, there are blowup sequences which converge to a non-constant time-independent (axisymmetric) harmonic map.
\end{abstract}

\tableofcontents

\section{Introduction}

The dynamic theory of nematic liquid crystals  was first proposed by Ericksen \cite{ericksen62} and Leslie \cite{leslie68}  in the 1960s, (see, e.g. \cite{Les, lin89}).  
More precisely, if the orientation order parameters of molecules are treated as a unit vector ${\d(x)}$, the classic static theory relies on the following Oseen-Frank free energy density
$$
2W(d,\nabla \d)
=k_1(\di \d)^2+k_2(\d\cdot\curl \d)^2+k_3|\d\times\curl \d|^2+(k_2+k_4)[\mbox{tr}(\nabla \d)^2-(\di \d)^2 ],
$$
where $k_i$, $i=1,\cdots, 3$ are the positive constants representing splay, twist, and bend affects respectively, with
$
k_2\geq |k_4|, \quad 2k_1\geq k_2+k_4.
$
The full Ericksen-Leslie system is given as follows

\begin{equation}\label{fwels}
\begin{cases}
\dot \u+\nabla P=\nabla\cdot\sigma-\nabla\cdot\left(\frac{\partial W}{\partial\nabla \d}\otimes\nabla \d\right),
\\
\nabla\cdot \u=0,\\
\nu \ddot{\d}=\gamma \d- g-\frac{\partial W}{\partial  \d}+\nabla\cdot\left(\frac{\partial W}{\partial\nabla \d}\right), \\
|\d|=1,
\end{cases}
\end{equation}
where $\u$ is the velocity field of underlying incompressible fluid,
$\dot{f}=f_t+\u\cdot \nabla f$ is the material derivative,
$\gamma$ is Lagrangian multiplier of the constraint $|\d|=1$, and the constant $\nu>0$ (we take $\nu=1$ in our paper).
Let
$$
D= \frac12(\nabla \u+(\nabla \u)^T),\quad \omega= \frac12(\nabla \u-(\nabla \u)^T)=\frac12(\partial_j\u^i-\partial_i\u^j),\quad N=\dot \d-\omega \d,
$$
represent the rate of strain tensor, skew-symmetric part of the strain rate, and the
rigid rotation part of director changing rate by fluid vorticity, respectively. The kinematic transport $g$ is given by
$$
g=\lambda_1 N +\lambda_2Dd
$$
which represents the effect of the macroscopic flow field on the microscopic structure. The material coefficients $\lambda_1$ and $\lambda_2$ reflect the molecular shape and the slippery part between fluid and particles. The first term of $g$ represents the rigid rotation of molecules, while the second term stands for the stretching of molecules by the flow.
The viscous (Leslie) stress tensor $\sigma$ has the following form
\beq\label{defsigma}
\sigma= \mu_1 (\d^TD\d)\d\otimes \d + \mu_2N\otimes \d   + \mu_3\d\otimes N+ \mu_4D + \mu_5(D\d)\otimes \d+ \mu_6\d\otimes (D\d) .
\eeq
These coefficients $\mu_i (1 \leq i \leq 6)$, depending on material and temperature, are called Leslie coefficients, and are related to certain local correlations in the fluid. The following relations are assumed in the literature.
\beq\label{parodis}
\lambda_1 =\mu_3-\mu_2,\quad \lambda_2 =\mu_6 -\mu_5,\quad \mu_2+ \mu_3 =\mu_6-\mu_5.
\eeq
The first two relations are compatibility conditions, while the third relation is called Parodi's relation, derived from Onsager reciprocal relations expressing the equality of certain relations between flows and forces in thermodynamic systems out of equilibrium. They also satisfy the following empirical relations (p.13, \cite{Les}) 
\begin{align}\label{alphas}
&\mu_4>0,\quad 2\mu_1+3\mu_4+2\mu_5+2\mu_6>0,\quad \lambda_1=\mu_3-\mu_2>0,\\
&  2\mu_4+\mu_5+\mu_6>0,\quad 4\lambda_1(2\mu_4+\mu_5+\mu_6)>(\mu_2+\mu_3+\gamma_2)^2\notag.
\end{align}
Note that the 4th relation  is implied by the 3rd together with the last relation.


 In this paper, we study the following Dirichlet energy density as the special case of Oseen-Frank energy density
$$
2W(\d,\nabla \d)=|\nabla \d|^2,
$$
and the system \eqref{fwels} becomes
\begin{equation}\label{fels}
\begin{cases}
 \u_t+\u\cdot \nabla \u+\nabla P=\nabla\cdot\sigma-\nabla\cdot\left(\nabla \d\odot\nabla \d\right),
\\
\nabla\cdot \u=0,\\
\ddot{\d}+\lambda_1(\dot \d-\omega \d)+\lambda_2 D\d=\Delta \d+\big(|\nabla \d|^2-|\dot{\d}|^2+\lambda_2\d^TD\d\big)\d, \\
|\d|=1.
\end{cases}
\end{equation}

In \cite{cai-wang20, HJLZ21, jiangluo17, Wang24} the global wellposedness of the classical solutions to the system \eqref{fels} has been proved under the smallness assumption  of the initial data for various choices of coefficients, using the dispersive, {dissipative} or damping mechanism.  However, for the general Cauchy problem, due to the wave structure of the director field, the large solution may develop finite time blowup \cite{CHX23}. The classical solutions cannot continue after the singularity formation, so it is meaningful to study the weak solutions.


We will consider the axisymmetric Poiseuille laminar flow via a tube for \eqref{fels}. More precisely, for $(x,y)\in \mathbb R^2$ and $r=\sqrt{x^2+y^2}$, we consider
\beq\label{poisup}
\u(x,y,t)=(0,0,v(r,t))^T,\quad P(x,y,t)=P(r,t),
\eeq
\beq\label {axsimd}
\d(x,y,t)=\big(\sin\phi(r,t) \cos \theta,\sin\phi(r,t)\sin \theta, \cos\phi(r,t) \big)^T.
\eeq
For simplicity, in this paper, we only consider a simple case of the coefficients in \eqref{fels}  
\beq\label{specialmus}
 \mu_1=\mu_5=\mu_6=0,\quad \mu_2=-1,\quad \mu_3=1,\quad \mu_4=1, \quad \lambda_1=\mu_3-\mu_2=2.
\eeq
By the Onsager-Parodi relation,
\[
\lambda_2=\mu_6-\mu_5=\mu_2+\mu_3=0.
\]
Then the system \eqref{fels} becomes the following system (see Appendix in \cite{CHX23} for details)
\begin{equation}\label{wels}
\begin{cases}
\displaystyle v_t=\frac{1}{r}\Big(rv_r+r\phi_t\Big)_r,\\
\\
 \displaystyle \phi_{tt}+2\phi_t=\frac{1}{r}\big(r\phi_r\big)_r-\frac{\sin(2\phi)}{2r^2}-v_r,
\end{cases}
\end{equation}
where $(r,t)\in (0,\infty)\times(0,\infty)$. The local wellposedness of classical solutions has been covered by \cite{jiangluo17} as a special case. When the initial data is small in Sobolev space $H^s, (s>2)$, the global solution has been constructed in \cite{jiangluo17}. 

In this paper, we study the existence and partial regularity of global weak solutions to the system \eqref{wels} with the initial and boundary conditions
\beq\label{welsr}
v(r,0)=v_0(r)\in H^1,\ \phi(r,0)=\phi_0(r)\in H^1,\ \phi_t(r,0)=\phi_1(r)\in L^2,\ 
v(0,t)=\phi(0,t)=0. 
\eeq
Here, and in the following, with slightly abusing notations, we use $H^s$ to denote the Sobolev space $H^s(\R^2)$ in 2D, which in the case of axisymmetric case, means
\[
\left\{f(r)\Bigg| \int_0^\infty \left(|f(r)|^2+\left|\frac{d^sf}{dr^s}(r)\right|^2 \right) rdr<\infty.\right\}.
\]
We use $\dot H^s$ for the homogeneous Sobolev space. 

One notices that the second equation of \eqref{wels} includes a damped
$O(3)$-$\sigma$ model. If we omit the coupling term $v_r$ and the damping term $\phi_t$, it becomes the axisymmetric $O(3)$-$\sigma$ model (when $k=1$) 
\beq\label{o3si}
\phi_{tt}=\frac{1}{r}\big(r\phi_r\big)_r-k^2\frac{\sin(2\phi)}{2r^2},\qquad k\geq 1,
\eeq
which can be derived from the wave map equations under a symmetry in dimension two. See \cite{Shatah01} for the model and existence of finite energy ($H^1$-type) solution.

The $O(3)$-$\sigma$ model from  wave map equations has been intensively studied.  
For smooth initial data $(\phi_0,\phi_1)$, existence of global regular solutions to \eqref{o3si} has been established in the pioneering works by Shatah-Tahvildar-Zadeh \cite{Shatah02} by assuming the initial energy $E(\phi)(0)$ is sufficiently small (see also \cite{ssb98} Theorem 8.1), where 
\[
\tilde E(\phi)(t)=\pi\int_{\mathbb R^+}\left(
(\phi_t)^2+(\phi_r)^2+\frac{k^2}{r^2}\sin^2 \phi\right)\, rdr.
\]
In \cite{struwe03}, Struwe showed that the singularities of equation \eqref{o3si} can only be developed by concentrating energy at the tip of a light cone by bubbling off at least one non-trivial harmonic map with finite energy. 
When $\phi(0,t)=0,\ \phi(\infty, t)=\pi,$
there are several approaches establishing examples of singularity formation of \eqref{o3si} with $\tilde{E}(\phi)(0)=4k \pi+\varepsilon$ by Rodnianski and Sterbenz \cite{RS} when $k\geq 4$, Krieger-Schlag-Tataru \cite{CST} when $k=1$, and \cite{RR} for all other cases. 
For solutions satisfying $\phi(0,t)=0,\ \phi(\infty, t)=0,$
 the global existence and scattering have been proved by a sequence of seminal works: C\^{o}te-Kenig-Merle \cite{CKM},
C\^{o}te-Kenig-Lawrie-Schlag \cite{CKL},  Jendrej-Lawrie \cite{JL}, when $$\tilde{E}(\phi)(0)\leq 8k \pi.$$
There are many other fundamental works on the existence and behavior of solutions for \eqref{o3si}, such as \cite{SU}.

In this paper, we study the properties of weak solutions and singularities for the nematic liquid crystals model \eqref{wels}. Our method depends on some existing methods for the $O(3)$-$\sigma$ model discussed in the last paragraph, and our recent progress in \cite{CHL19, CHX23}.  When one considers more general Oseen-Frank energy in \eqref{fwels}, the formation of cusp singularity and global weak solutions have been constructed in \cite{CHL19} in dimension one by utilizing the Poiseuille flow. Here the cusp singularity forms due to the quasi-linearity in the ${\bf d}$ or $\phi$ equation produced by the general Oseen-Frank energy. For \eqref{wels} with a semilinear wave equation on $\phi$, one needs to cope with a totally different singularity comparing to the one in \cite{CHL19} at the origin due to the geometric effect in multiple space dimension, to some extent similar to the $O(3)$-$\sigma$ model. In \cite{CHX23}, a singularity formation of \eqref{wels} has been constructed following the frame work of \cite{RS}  for $k\geq 4$ (see \eqref{o3si} for definition of $k$).

The formation of this kind of singularity in \eqref{wels} is the major difficulty in studying the global existence of finite energy solution, our first major result. We use the Ginzburg-Landau approximation and fixed-point arguments to establish the global existence of weak solutions, to go beyond the singularity.

\medskip
Our first main result concerns the existence of global weak solutions  to the system \eqref{wels}.
\begin{definition}[weak solution]\label{wsdef1}
$v,\phi$ is called a weak solution of \eqref{wels} if $v\in L^\infty((0,t), L^2)\cap L^2((0,t), \dot H^1),\quad \phi\in L^\infty((0,t), \dot H^1),\quad  \phi_t\in L^\infty((0,t), L^2)$ satisfy for any compactly supported test function $\psi_1\in X_1:=\{f:\R_+\times(0,T)\to \R,  \mid f\in L^\infty((0,T),H^1),~f_t\in L^\infty((0,T),L^2)\}$, $\psi_2\in X_2:=\{f\in C^\infty(\R^+\times\R^+)\mid \exists M>0, \forall t\geq 0, f(r,t)=0$, when $r>M\}$.
\begin{equation}\label{eqnphi}
\int_0^\infty\int_0^\infty\left( -\phi_t\psi_{1,t}+2\phi_t\psi_1+\phi_r\psi_{1,r}+v_r\psi_1+\frac{\sin(2\phi)}{2r^2}\psi_1\right) rdr=0,
\end{equation}
and
\begin{equation}\label{eqnv}
\int_0^\infty\int_0^\infty -v\psi_{2,t}+(v_r+\phi_t)\psi_{2,r} rdrdt=0.
\end{equation}
\end{definition}

\btm\label{thm1}
For any $t>0$, the system \eqref{wels} with initial and boundary conditions \eqref{welsr} has a weak solution 
\beq\notag
v\in L^\infty((0,t), L^2)\cap L^2((0,t), H^1),\quad \phi\in L^\infty((0,t), H^1),\quad  \phi_t\in L^\infty((0,t), L^2).
\eeq
\etm


The second main theorem studies the partial regularity of the solutions.  
Inspired by \cite{CHX23} and \cite{CHL19}, we introduce a quantity 
$$h(r,t)=\frac{1}{r}\int_0^rv(R,t)\,RdR.$$
By the definition of $h$ and \eqref{welsr}, it holds
\beq\label{welshr}
h(0,t)=0,\quad h(r,0)=\frac{1}{r}\int^r_0v_0(R)\,RdR\triangleq h_0(r).
\eeq
Direct computation implies
$$
\left(\int_0^rv(R,t)\,R\, dR\right)_t=rv_r+r\phi_t, \quad v_r=\left(\frac{(rh)_r}{r}\right)_r
=h_{rr}+\frac{h_r}{r}-\frac{h}{r^2},
$$
which implies
\beq\label{eqnvth}
h_t=v_r+\phi_t=\frac1r(rh_r)_r-\frac{h}{r^2}+\phi_t.
\eeq
Hence the system \eqref{wels} becomes
\begin{equation}\label{welss}
\begin{cases}
\displaystyle h_t=\frac{1}{r}(rh_r)_r-\frac{h}{r^2}+\phi_t,\\
\\
\displaystyle \phi_{tt}+\phi_t=\frac{1}{r}(r\phi_r)_r-\frac{\sin(2\phi)}{2r^2}-h_t.
\end{cases}
\end{equation}
Notice that we replaced the $v_r$ by $h_t-\phi_t$ in the second equation. It is all because $h_t$ has much better ``regularity" than $v_r$, which is very crucial in our results on global existence and partial regularity of solution. This idea has been utilized in \cite{CHL19, CHX23}. 

We define the local energy density and local energy
$$
 e(r,t)=\frac{1}{2}|\phi_r|^2+\frac{1}{2}|\phi_t|^2+\frac{\sin^2\phi}{2r^2},\quad E(R,t)=\int^R_0e(r,t)\,rdr.
$$
For any $T>0$, suppose $R> 0$, $s>0$ be any constants satisfying $T=R+s$. We consider the following region 
$$
K^T(s,\tau)= \{0\leq s-\tau\leq t\leq s,\ 0\leq r\leq T-t\},\quad 
M^T(s,\tau)= \{0\leq s-\tau\leq t\leq s,\ r=T-t\}
$$
for any $0<\tau\leq s< T$.
Denote
$$
\text{Flux}(s,s-\tau)= \int_{s-\tau}^s\big[e(T-t,t)-\phi_r\phi_t(T-t,t)\big]\,(T-t)dt.
$$
It is easy to see that Flux$(s,s-\tau)\ge 0$.

For simplicity, we introduce our results on \eqref{welss}, which also work for \eqref{wels}. 

\begin{theorem}\label{thma}
For any smooth local solutions $(h,\phi)$ to the initial-boundary value problem \eqref{welss} and \eqref{welsr}, if there exist $T_0>0$ and $\delta>0$, such that when $0<T_0-\delta<t<T_0$ and $R>T_0-t$, it holds 
\beq\label{ineq1}
E(R,t)<\epsilon_0,
\eeq
for some given $\epsilon_0>0$, then $(h,\phi)$ can be extended continuously over $T_0$.
\end{theorem}

{
\begin{remark}
(1) Since the system is a special case of the full system studied in \cite{jiangluo17}, the local well-posedness of classical solutions follow directly from the theorem 1.1 of \cite{jiangluo17}.\\
\smallskip
(2) Due to technical difficulties, the higher regularity of the solutions is still open. For axisymmetric $O(3)$-$\sigma$ model \eqref{o3si}, it has been proved by Shatah-Tahvildar-Zadeh in \cite{Shatah02} with small energy.
\end{remark}
}


Another interesting result is on the structure of the first possible blowup at finite time. 

\btm\label{thmb}
Suppose that the solution $(\phi,h)$ to the system \eqref{welss} blows up at $T_0>0$. Then there exist  sequences $R_i\rightarrow 0^+$, $T_i\rightarrow T_0^-$ and the scaling
$$
\phi_i(r,t)=\phi(R_ir,T_i+R_it),\quad h_i(r,t)=h(R_ir,T_i+R_it)
$$
such that 
$$
\phi_i(r,t)\rightarrow\phi_\infty(r,t)\text{ in }H^1_{loc}((0,\infty)\times(-1,1)),\quad 
h_i(r,t)\xrightarrow{w}0\text{ in }H^1_{loc}$$
where $\phi_\infty(r,t)$ is a non-constant time independent axisymmetric harmonic map satisfying
$$(\phi_\infty)_{rr}+\frac{(\phi_\infty)_r}{r}-\frac{\sin(\phi_\infty)\cos(\phi_\infty)}{r^2}=0.$$
\etm

One should notice that in the blowup sequence, we choose the scaling of wave in time variable instead of scaling of parabolic equation since the possible singularities of system are generated by the nonlinearity of the wave equation according to our paper \cite{CHX23}. 


{
Finally, when the initial energy less than the energy of nontrivial harmonic map, we show that the system \eqref{welss} has no $C^0$ singularities as in \cite{CHX23}.
\begin{theorem}\label{extthm}
Suppose that initial energy of the system \eqref{welss} satisfies the following condition
\beq\label{sminen}
E_0=\int\left(|(\phi_0)_r|^2+|\phi_1|^2+\frac{\sin^2(\phi_0)}{r^2}+|(h_0)_r|^2+\frac{h_0^2}{r^2}\right)\,rdr<4.
\eeq
Then the weak solutions $(h,\phi)$ to the initial-boundary value problem \eqref{welss}, \eqref{welsr} and \eqref{welshr} together with another extra boundary condition $\phi(\infty, t)=0$ have no $C^0$ blowup at any time $0<t<\infty$.
\end{theorem}
}

The remaining of the paper is organized as follows. In Section 2, we construct global weak solutions by the Ginzburg-Landau approximation and the fixed-point argument. In Section 3, we obtain the global and local energy estimates which are necessary in proving partial regularity. In Section 4, we prove the partial regularity theorem \ref{thma}. In Section 5, we construct blowup sequences at the first singular time and prove Theorem \ref{thmb}. In Section 6, we prove Theorem \ref{extthm}.



\section{Construction of global weak solutions}

In this section, we will investigate the existence of global weak solutions to the system \eqref{wels} with initial and boundary conditions \eqref{welsr}. We first rewrite the the system back to the Cartesian coordinate system as follows
\begin{equation}\label{fels2}
\begin{cases}
 v_t-\frac12\Delta v=\nabla\cdot\mathcal T_1+\nabla\cdot\mathcal T_2,
\\
\d_{tt}+2 \d_t-2\omega \d=\Delta \d+\big(|\nabla \d|^2-|\d_t|^2\big)\d, \\
|\d|=1.
\end{cases}
\end{equation}
Where 
\beq\label{defdw}
\d=\left(d^1,d^2,d^3\right)^T,\quad 
2\omega=
\left(
\begin{array}{ccc}
0&0&-v_x\\
0&0&-v_y\\
v_x&v_y&0\\
\end{array}
\right),
\eeq
\beq\label{defcT1}
\mathcal T_1(\d,d_t)=\left(
d^3d^1_t-d^1d^3_t,\ d^3d^2_t-d^2d^3_t
\right)^T,
\eeq
\beq\label{defcT2}
\mathcal T_2(\d, \nabla v)=\frac12\big(
(|d^3|^2+|d^1|^2)v_x+d^1d^2v_y,\ (|d^3|^2+|d^2|^2)v_y+d^1d^2v_x
\big)^T.
\eeq
The last two terms come from the term $\d\otimes N- N\otimes \d$ in viscous stress tensor $\sigma$.
It is not hard to verify the first equation of system \eqref{fels2} is still uniform parabolic.

\begin{definition}[weak solution]\label{wsdef2}
$(v, \d)$ is a weak solution of \eqref{fels2} on $\R^2\times(0,T)$ if \\
$v\in L^\infty((0,T),L^2)\bigcap L^2((0,T),\dot H^1)$, $\d\in L^\infty((0,T),\dot H^1), \d_t\in L^\infty((0,T), L^2)$ satisfy for any compactly supported test function $\varphi_1\in Y_1:=\{f:\R^2\times(0,T)\to \R,  \mid f\in L^\infty((0,T),H^1(\R^2)),~f_t\in L^\infty((0,T),L^2(\R^2))\}$, $\varphi_2\in C_c^\infty(\R^2\times \R^+)$,
\begin{equation}\label{ws2}
\int_0^T\int_{\R^2}\left\{-\d_t\varphi_{1,t}+2\d_t\varphi_1-2\omega\d\varphi_1+\nabla\varphi_1\cdot\nabla\d-(|\nabla \d|^2-|\d_t|^2)\d\varphi_1\right\} d\x dt
=0,
\end{equation}
and
\begin{equation}\label{ws1}
\int_0^T\int_{\R^2} \left\{-v\varphi_{2,t}+\frac12\nabla v\cdot\nabla \varphi_2+(\mathcal T_1+\mathcal T_2)\cdot\nabla\varphi_2 \right\}d\x dt=0.
\end{equation}

\end{definition}

About the weak solutions of \eqref{wels} and \eqref{fels2}, we have the following proposition.
\begin{proposition}
Assume $(v,\d)$ is a weak solution of \eqref{wels} in the sense of Definition \ref{wsdef2} satisfying 
\begin{itemize}
\item The solution is axisymmetric, and
\begin{equation}\label{def_d}
\d(x,y,t)=\big(\sin\phi(r,t) \cos \theta,\sin\phi(r,t)\sin \theta, \cos\phi(r,t) \big)^T,
\end{equation} 
\item  $\phi|_{r=0}=0$. 
 \end{itemize}
Then $(v,\phi)$ is a weak solution of \eqref{fels2} in the sense of Definition \ref{wsdef1}.
 \end{proposition}
\pf
Assume $(v, \d)$ is a weak solution of \eqref{fels2}. By \eqref{def_d}, 
 $\d\in \dot H^1$ leads to $\phi\in \dot H^1$. Therefore,
$v\in L^\infty((0,t), L^2)\cap L^2((0,t), \dot H^1),\quad \phi\in L^\infty((0,t), \dot H^1),\quad  \phi_t\in L^\infty((0,t), L^2)$.

Considering \eqref{def_d}, the first component of \eqref{ws2} leads to 
\begin{align*}
\int_0^{2\pi}\iint_{\R_+\times\R_+}&\left\{\cos\phi\cos\theta[-\phi_t\varphi_t+2\phi_t\varphi+v_r\varphi+\phi_r\varphi_r]r-\sin\phi\cos\theta \left(\phi_r^2-\phi_t^2+\frac{\sin^2\phi}{r^2}\right)\varphi r\right.\\
&\left.+\varphi\cos\theta \frac{\sin\phi}{r^2}r\right\}dr dt d\theta
=0,
\end{align*}
for any test function $\varphi\in Y_1$. 
Use the axisymmetric assumption of $\phi$, define 
$$f(r,t)=\int_0^{2\pi}\cos\theta\varphi(r,\theta, t)d\theta.$$ Then $f\in X_1$. $f$ is a test function of Definition \ref{wsdef1}. Then we have
\begin{equation}
\iint_{\R_+\times\R_+}\left\{\cos\phi[-\phi_tf_t+2\phi_tf+v_rf+\phi_rf_r]r-\sin\phi \left(\phi_r^2-\phi_t^2\right)f r+ \frac{\cos^2\phi\sin\phi}{r^2}f r\right\}dr dt 
=0,
\end{equation}
which is written as
\begin{equation}\label{1stcomp}
\iint_{\R_+\times\R_+}\left\{\cos\phi\left[-\phi_tf_t+2\phi_tf+v_rf+\phi_rf_r+\frac{\sin(2\phi)}{2r^2}f\right]r-\sin\phi \left(\phi_r^2-\phi_t^2\right)f r\right\}dr dt 
=0.
\end{equation}
Since $\varphi$ is an arbitrary test function, $f$ is also arbitrary. 

The third component of \eqref{ws2} leads to 
\begin{equation}
\int_0^{2\pi}\iint_{\R_+\times\R_+}\left\{\sin\phi[\phi_t\varphi_t-2\phi_t\varphi-v_r\varphi-\phi_r\varphi_r]r-\cos\phi \left(\phi_r^2-\phi_t^2+\frac{\sin^2\phi}{r^2}\right)\varphi r\right\}dr dt d\theta
=0.
\end{equation}
Defining $f(r,t)=\int_0^{2\pi}\varphi(r,\theta, t)d\theta$ as the test function, we have,
\begin{equation}\label{2ndcomp}
\iint_{\R_+\times\R_+}\left\{\sin\phi\left[\phi_tf_t-2\phi_tf-v_rf-\phi_rf_r-\frac{\sin(2\phi)}{2r^2}f\right]r-\cos\phi \left(\phi_r^2-\phi_t^2\right)f r\right\}dr dt 
=0.
\end{equation}

Then we take test function $f$ to be $\cos\phi\psi$ and $\sin\phi\psi$ separately in \eqref{1stcomp} and \eqref{2ndcomp}, where $\psi\in X_1$ is axisymmetric and we obtain
\begin{align*}
&\iint_{\R_+\times\R_+}\left\{\cos^2\phi\left[-\phi_t\psi_t+2\phi_t\psi+v_r\psi+\phi_r\psi_r+\frac{\sin(2\phi)}{2r^2}\psi\right]r-\sin\phi\cos\phi \left(\phi_r^2-\phi_t^2\right)\psi r\right\}dr dt \\
&+\iint_{\R_+\times\R_+}\cos\phi\sin\phi(\phi_t^2-\phi_r^2)\psi rdrdt
=0,
\end{align*}
and
\begin{align*}
&\iint_{\R_+\times\R_+}\left\{\sin^2\phi\left[\phi_t\psi_t-2\phi_t\psi-v_r\psi-\phi_r\psi_r-\frac{\sin(2\phi)}{2r^2}\psi\right]r-\sin\phi\cos\phi \left(\phi_r^2-\phi_t^2\right)\psi r\right\}dr dt \\
&+\iint_{\R_+\times\R_+}\cos\phi\sin\phi(\phi_t^2-\phi_r^2)\psi rdrdt=0.
\end{align*}
Then by taking the difference of these two equations, we obtain
\[
\iint_{\R_+\times\R_+}\left[-\phi_t\psi_t+2\phi_t\psi+v_r\psi+\phi_r\psi_r+\frac{\sin(2\phi)}{2r^2}\psi\right]rdr dt=0.
\]
So $\phi$ satisfies the equation \eqref{eqnphi}. 

Next, we consider the weak solution for $v$. If $v$ satisfies \eqref{ws1}, then under the axisymmetric assumption, for any $\varphi\in Y_2$
\begin{align*}
\mathcal T_1+\mathcal T_2=
\left(
\begin{array}{c}
\cos\theta\phi_t+\frac12v_r\cos\theta\cr\cr
\sin\theta\phi_t+\frac12v_r\sin\theta
\end{array}
\right),
\end{align*}
\begin{align*}
\nabla\varphi=\left(
\begin{array}{c}
\varphi_r\cos\theta-\varphi_\theta\frac{\sin\theta}{r}\\
\varphi_r\sin\theta+\varphi_\theta\frac{\cos\theta}{r}
\end{array}
\right),
\end{align*}
thus
\[
(\mathcal T_1+\mathcal T_2)\cdot\nabla\varphi=\varphi_r\phi_t+\frac12v_r\varphi_r.
\]
\[
\nabla v\cdot\nabla \varphi=\varphi_rv_r.
\]
So we have
\[
\int_0^T\int_0^{2\pi}\int_{\R^+} (-v\varphi_t+\varphi_rv_r+\varphi_r\phi_t) rdrd \theta dt=0.
\]
Define $\psi(r,t)=\int_0^{2\pi}\varphi(r\cos\theta, r\sin\theta, t)d\theta$. Then $\psi\in X_2$, and $v$ satisfies
\[
\int_0^\infty\int_0^\infty \left[-v\psi_{t}+(v_r+\phi_t)\psi_{r} r\right]drdt=0.
\]
Therefore, $(v,\phi)$ is a weak solution of \eqref{wels} in the sense of Definition \ref{wsdef1}.

\endpf

We can derive the following energy estimate.

\blm\label{lemma2.1}
Suppose that $v$, $\d$ are smooth solutions to the system \eqref{fels2}. Then it holds the energy equality
\beq\label{engestw1}
\frac{d}{dt} \frac12\int \left(|v|^2+|\d_t|^2+|\nabla \d|^2\right)\,d\x 
+\frac12\int|\nabla v|^2\,d\x=-2\int|\d_t-\omega \d|^2\,d\x,
\eeq
where $\x=(x,y)\in \R^2$, $\int\,d\x=\int_{\R^2}\,d\x$ for simplicity.
\elm

\pf Multiplying the first equation of \eqref{fels2} by $v$ and integrating by parts over $\R^2$, we have 
\beq\label{lem1pf1}
\frac{d}{dt} \frac12\int |v|^2\,d\x +\frac12\int |\nabla v|^2\,d\x
=-\int \mathcal T_1(\d,\d_t)\cdot \nabla v\,d\x-\int \mathcal T_2(\d,\nabla v)\cdot \nabla v\,d\x.
\eeq
Direct computation implies
\beq\label{lem1pf2}
-\int \mathcal T_1(\d,\d_t)\cdot \nabla v\,d\x
=-\int\left[ \big(d^3d^1_t-d^1d^3_t\big)v_x+\big(d^3d^2_t-d^2d^3_t\big)v_y\right]\,d\x,
\eeq
\beq\label{lem1pf3}
-\int \mathcal T_2(\d,\nabla v)\cdot \nabla v\,d\x
=-\frac12\int\left[ \big(|d^3|^2+|d^1|^2\big)|v_x|^2+\big(|d^3|^2+|d^2|^2\big)|v_y|^2+2d^1d^2v_xv_y\right]\,d\x.
\eeq
Multiplying the second equation of \eqref{fels2} by $\d_t$ and integrating by parts over $\R^2$, we have 
\beq\label{lem1pf4}
\frac{d}{dt} \frac12\int \left(|\d_t|^2+|\nabla \d|^2\right)\,d\x +2\int |\d_t|^2\,d\x
=2\int\omega \d\cdot \d_t\,d\x.
\eeq
Direct computation implies
\beq\label{lem1pf5}
2\int\omega \d\cdot \d_t\,d\x
=-\int\left[ \big(d^3d^1_t-d^1d^3_t\big)v_x+\big(d^3d^2_t-d^2d^3_t\big)v_y\right]\,d\x.
\eeq
Combining \eqref{lem1pf1}-\eqref{lem1pf5}, we obtain
\beq\notag
\begin{split}
&\frac{d}{dt} \frac12\int \left(|v|^2+|\d_t|^2+|\nabla \d|^2\right)\,d\x 
+\int\left(\frac12|\nabla v|^2+2 |\d_t|^2\right)\,d\x\\
=&-2\int\left[ \big(d^3d^1_t-d^1d^3_t\big)v_x+\big(d^3d^2_t-d^2d^3_t\big)v_y\right]\,d\x\\
&-\frac12\int\left[ \big(|d^3|^2+|d^1|^2\big)|v_x|^2+\big(|d^3|^2+|d^2|^2\big)|v_y|^2+2d^1d^2v_xv_y\right]\,d\x\\
=&-2\int\left[\left(d^1_t+\frac12d^3v_x\right)^2+\left(d^2_t+\frac12d^3v_y\right)^2
+\left(d^3_t-d^1v_x-d^2v_y\right)^2\right]\,d\x+2\int |\d_t|^2\,d\x,
\end{split}
\eeq
from which we can derive the energy inequality \eqref{engestw1}. This completes the proof of the lemma.
\endpf

\medskip
To construct the global weak solution, we use the following Ginzburg-Landau approximation of system \eqref{fels2}, which releases the constrain $|\d|=1$ by the energy $\frac{1}{4\v^2}(|\d^\v|^2-1)^2$ for $\v>0$
\begin{equation}\label{fels3}
\begin{cases}
 v^\v_t-\frac12\Delta v^\v=\nabla\cdot\mathcal T_1(\d^\v,\d^\v_t)+\nabla\cdot\mathcal T_2(\d^\v,\nabla v^\v),
\\
\d^\v_{tt}+2 \d^\v_t-2\omega^\v \d^\v=\Delta \d^\v-\frac{|\d^\v|^2-1}{\v^2}\d^\v,
\end{cases}
\end{equation} 
where  $\omega^\v$ are the same form as in \eqref{defdw} by replacing $v$ by $v^\v$. Repeat the procedure in proof of Lemma \ref{lemma2.1}, we obtain the following energy estimate of system \eqref{fels3}.

\blm\label{lemma2.2}
Suppose that $v^\v$, $\d^\v$ are smooth solutions to the system \eqref{fels3}. Then it holds
\beq\label{engestw2}
\begin{split}
&\frac{d}{dt} \frac12\int \left(|v^\v|^2+|\d^\v_t|^2+|\nabla \d^\v|^2+\frac{1}{4\v^2}(|\d^\v|^2-1)^2\right)\,d\x\\
=&-\frac12\int|\nabla v^\v|^2\,d\x-2\int|\d^\v_t-\omega^\v \d^\v|^2\,d\x.
\end{split}
\eeq
\elm

We first construct the global solutions to the system \eqref{fels3} with initial conditions
\beq\label{initial2}
v^\v(x,y,0)=v_0(r),\quad \d^\v(x,y,0)=\d_0(x,y)=\big(\sin\phi_0(r)\cos\theta,\sin\phi_0(r)\sin\theta,\cos\phi_0(r)\big)^T,
\eeq
\beq\label{initial3}
\d^\v_t(x,y,0)=\d_1(x,y)=\big(\cos\phi_0(r)\cos\theta,\cos\phi_0(r)\sin\theta,-\sin\phi_0(r)\big)^T\phi_1(r).
\eeq
Since $\phi_0\in H^1$, it is reasonable to assume $\lim\limits_{r\rightarrow \infty}\d^\v_0(r)=(0,0,1)^T:=\e_1$
and
\beq\label{initial4}
v_0\in H^1,\quad \d_0-\e_1\in H^1,\quad \d_1\in L^2,\quad 
\lim\limits_{r\rightarrow \infty}(v,\d)(r,t)=(0,\e_1).
\eeq
Denote
\beq\label{defcX}
\begin{split}
\mathcal X:=& \big\{\mathbf{f}(r)-\e_1\in L^\infty((0,t), H^2(\R^2,\R^3)),\  \mathbf{f}_t(r,t)\in L^\infty((0,t), H^1(\R^2,\R^3)),\\
& \mathbf{f} \mbox{ is axisymmetric without swirl}\big\},
\end{split}
\eeq
with norm
\beq\label{defnmX}
 \|\mathbf{f}\|_{\mathcal X}:=\|\mathbf{f}-\e_1\|_{ L^\infty((0,t), H^2(\R^2,\R^3))}+\|\mathbf{f}_t\|_{L^\infty((0,t), H^1(\R^2,\R^3))}
\eeq
with distance 
\beq
\rho_{\mathcal X}(\mathbf f_1,\mathbf f_2):= \sup\limits_{0\leq s\leq t}
\left(\|\mathbf f_1-\mathbf f_2\|_{H^1}+\|(\mathbf f_1)_t-(\mathbf f_2)_t\|_{2}\right).
\eeq
Here we denote $\|\cdot\|_k$ as the $L^k$-norm for $\R^2$.
For any $K>0$, denote 
\beq\label{defcXK}
\mathcal X_K:= \big\{\mathbf f\in \mathcal X,\ \mathbf f(x,y,0)=\d_0,~~\mathbf f_t(x,y,0)=\d_1, ~~ \|\mathbf f\|_{\mathcal X}<K\big\}.
\eeq

\begin{proposition}\label{prop2.3}
For any $t>0$, the Cauchy problem of \eqref{fels3} with initial and boundary conditions \eqref{initial2}-\eqref{initial4} has a weak solution 
\beq\notag
v^\v\in L^\infty((0,t), L^2(\R^2,\R))\cap L^2((0,t), H^1(\R^2,\R)),
\eeq
\beq\notag
\d^\v\in L^\infty((0,t), H^1(\R^2,\R^3)),\quad \d^\v_t\in L^\infty((0,t), L^2(\R^2,\R^3)),
\eeq 
for any fixed $\v>0$. In additionally, the solutions $(v^\v,\d^\v)$ depend only on $(r,t)$ and $\d^\v$ is axisymmetric without swirl.
\end{proposition}

\pf Without loss of generality, we assume that $\v=1$ and omit the upper index $\v$ in the proof. 
Let's first assume that the initial data are little better as follows
\beq\label{initial5}
v_0\in H^1(\R^2,\R),\quad \d_0-\e_1\in H^2(\R^2,\R^3),\quad \d_1\in H^1(\R^2,\R^3).
\eeq
to solve the following approximate system for any $\delta >0$
\begin{equation}\label{felsd}
\begin{cases}
 v_t-\frac12\Delta v-\nabla\cdot\mathcal T_2(\d,\nabla v)=\nabla\cdot\mathcal T_1(\d,\d_t),
\\
\d_{tt}-\delta\Delta \d_t+2 \d_t-\Delta \d=2\omega(v) \d-(|\d|^2-1)\d,
\end{cases}
\end{equation} 
where  $\omega(v)$ is the same form as in \eqref{defdw}.

We solve \eqref{felsd} by fixed-point argument. More precisely, we consider the following linear system
\begin{equation}\label{fels4}
\begin{cases}
 v_t-\frac12\Delta v-\nabla\cdot\mathcal T_2(q,\nabla v)=\nabla\cdot\mathcal T_1(q,q_t),
\\
\d_{tt}-\delta\Delta \d_t+2 \d_t-\Delta \d=2\omega(v) q-(|q|^2-1)q,
\end{cases}
\end{equation} 
with the regularized initial data $v|_{t=0}=J^\delta v_0, \d|_{t=0}=J^\delta\d_0$, where $J^\delta$ is the smoothing operator defined as $J^\delta f(r)=\int_{\R^+} f(s) G_\delta(r-s)  sds$, $G$ is a smooth function supported on a unit disk satisfying $\int_{\R^+}G(r) rdr=1$, and $G_\delta(r)=\frac{1}{\delta^2}G(\frac{r}{\delta})$.

For any $q\in \mathcal X_K$, we solve $v$ from the first equation and then solve $\d$ from the second equation. By the standard theory of parabolic equations, the following map is well-defined 
\beq\notag
\d:=\mathcal L(q). 
\eeq
It is not hard to see that the solutions $(v,\d)$ only depend on $(r,t)$ and $\d$ is axisymmetric without swirl. 
The proof will be divided into several claims.

\medskip
\noindent{\it Claim 1. $v\in L^{\infty}((0,t),H^1)\cap L^{2}((0,t),H^2)$.} 
\medskip

\noindent Indeed, multiplying the first equation of \eqref{fels4} by $v$ and integrating by parts over $\R^2$, we have 
\beq\label{prop1pf00}
\frac{d}{dt} \frac12\int |v|^2\,d\x +\frac12\int |\nabla v|^2\,d\x+\int \mathcal T_2(q,\nabla v)\cdot\nabla v\,d\x
=-\int  \mathcal T_1(q,q_t)\cdot\nabla  v\,d\x.
\eeq
Similar to the computation of \eqref{lem1pf3}, we obtain
\beq\notag
\int \mathcal T_2(q,\nabla v)\cdot\nabla v\,d\x\geq 0.
\eeq
The term on the right side can be estimated as follows
\beq\notag
\begin{split}
\left|\int  \mathcal T_1(q,q_t)\cdot\nabla  v\,d\x\right|
\leq C\int |q||q_t||\nabla  v|\,d\x
\leq \frac14\int|\nabla v|^2\,d\x+C\|q\|_{\infty}^2\|q_t\|_2^2.
\end{split}
\eeq
The constant $C>0$ will be changed line by line.
Putting both estimates into \eqref{prop1pf00} and integrating over $(0,t)$, we obtain
\beq\label{prop1pf01}
\begin{split}
&\sup\limits_{0\leq s\leq t}\int |v|^2\,d\x +\int_0^t\int |\nabla v|^2\,d\x\\
\leq& \|v_0\|_2^2+C\int_0^t\|q\|_{\infty}^2\|q_t\|_2^2\,ds\\
\leq& \|v_0\|_2^2+C\int_0^t\|q\|_{H^2}^2\|q_t\|_2^2\,ds\\
\leq& \|v_0\|_2^2+CtK^4,
\end{split}
\eeq
where we have used the Sobolev inequality in the second inequality.

Since $q\in \mathcal{X}_K$, we have $q\in L^\infty ((0,t),C^{\alpha})$ for some $0<\alpha<1$. By the standard theory of general parabolic equation, it holds
\beq\label{prop1pf02}
\|v\|_{L^2 ((0,t),H^2)}^2
\leq C\|v\|_{L^2 ((0,t),L^2)}^2+C\|v_0\|_{2}^2+C\|\nabla\cdot\mathcal T_1(q,q_t)\|_{L^2 ((0,t),L^2)}^2.
\eeq
Direct computation implies
\beq\label{prop1pf03}
\begin{split}
&\int_0^t\int | \nabla \cdot \mathcal T_1(q,q_t)|^2\,d\x ds\\
\leq &\int_0^t\int \left(|\nabla q|^2|q_t|^2+| q|^2|\nabla q_t|^2\right)\,d\x ds\\
\leq & C\int_0^t\left(\|\nabla q\|_4^4+\|q_t\|_4^4+\|q\|_{\infty}^2\|\nabla q_t\|_2^2 \right)\,ds\\
\leq &C\int_0^t\left(\|\nabla q\|_2^2\|\Delta q\|_2^2+\|q_t\|_2^2\|\nabla q_t\|_2^2+\| q\|_{H^2}^2\|\nabla q_t\|_2^2 \right)\,ds\\
\leq &CtK^4,
\end{split}
\eeq
where we have used the { Ladyzhenskaya inequality 
$$\|\nabla q\|_4^4\leq C\|\nabla q\|_2^2\|\Delta q\|_2^2$$
 }and Sobolev inequality in the forth inequality. This combining with \eqref{prop1pf01} and \eqref{prop1pf02} implies 
\beq\label{prop1pf0}
\|v\|_{L^2 ((0,t),H^2)}^2
\leq C\|v_0\|_{2}^2+CtK^4.
\eeq

Multiplying the first equation of \eqref{fels4} by $\Delta v$ and integrating by parts over $\R^2$, we have 
\beq\label{prop1pf1}
\frac{d}{dt} \frac12\int |\nabla v|^2\,d\x +\frac12\int |\Delta v|^2\,d\x=
-\int \nabla \cdot \mathcal T_1(q,q_t)\Delta  v\,d\x-\int \nabla\cdot\mathcal T_2(q,\nabla v)\Delta v\,d\x.
\eeq
By the estimate \eqref{prop1pf02}, we have
\beq\label{prop1pf2}
\begin{split}
\left|\int \nabla \cdot \mathcal T_1(q,q_t)\Delta  v\,d\x\right|
\leq \frac18\int|\Delta v|^2\,d\x+C\int | \nabla \cdot \mathcal T_1(q,q_t)|^2\,d\x 
\leq \frac18\int|\Delta v|^2\,d\x+CK^4.
\end{split}
\eeq
For the second term on the right side of \eqref{prop1pf1}, it holds 
\beq\label{prop1pf3}
\begin{split}
&\left|\int \nabla\cdot\mathcal T_2(q,\nabla v)\Delta v\,d\x\right|\\
\leq &C\int\left(|q||\nabla q||\nabla v|+ |q|^2|\nabla^2 v|\right)|\Delta v|\,d\x\\
\leq& \frac18\int|\Delta v|^2\,d\x+C\left[\| q\|_{\infty}^2\|\nabla q\|_4^2\|\nabla v\|_4^2+\|q\|_{\infty}^4\|\nabla^2  v\|_2^2 \right]\\
\leq &\frac18\int|\Delta v|^2\,d\x+C\left[\| q\|_{H^2}^2\|\nabla q\|_2\|\Delta q\|_{2}\|\nabla v\|_2\|\Delta v\|_2+\|q\|_{H^2}^4\|\Delta  v\|_2^2 \right]\\
\leq &\frac14\int|\Delta v|^2\,d\x+CK^{8}\|\nabla v\|_2^2
+CK^4\|\Delta v\|_2^2.
\end{split}
\eeq
Combining \eqref{prop1pf1}-\eqref{prop1pf3} and integrating over $(0,t)$ with respect to time, we have
\beq\label{prop1pf31}
\begin{split}
&\sup\limits_{0\leq s\leq t}\int |\nabla v|^2(s)\,d\x +\int_0^t\int |\Delta v|^2\,d\x ds\\
\leq& \|\nabla v_0\|_2^2+CK^4t+CK^8\int_0^t\| v\|_2^2\,ds+CK^4\int_0^t\|\Delta v\|_2^2\,ds\\
\leq & \|\nabla v_0\|_2^2+CK^4t+C(K^8+K^4)(\|v_0\|_{2}^2+tK^4),
\end{split}
\eeq
where we have used the \eqref{prop1pf0} in last step. Combining \eqref{prop1pf31} with \eqref{prop1pf01} and \eqref{prop1pf0}, we prove the claim. 

\medskip
\noindent{\it Claim 2. $\d-\e_1 \in \mathcal{X}_K$.} 
\medskip

\noindent Indeed, operating $\nabla$ on the second equation of \eqref{fels4}, multiplying both sides by $\nabla \d_t$ and integrating by parts over $\R^2$, we have 
\beq\label{prop1pf4}
\begin{split}
&\frac{d}{dt} \frac12\int \left(|\nabla \d_t|^2+|\Delta \d|^2\right)\,d\x +2\int |\nabla\d_t|^2\,d\x+\delta\int|\Delta \d_t|^2\,d\x\\
=&\int 2\nabla(\omega(v) q)\cdot \nabla\d_t\,d\x-\int\nabla\left((|q|^2-1)q\right)\cdot \nabla\d_t\,d\x.
\end{split}
\eeq
Direct computation implies
\beq\label{prop1pf5}
\begin{split}
&\left|\int 2\nabla(\omega(v) q)\cdot \nabla\d_t\,d\x-\int\nabla\left((|q|^2-1)q\right)\cdot \nabla\d_t\,d\x\right|\\
\leq& C\int\left(|\nabla^2 v||q||+|\nabla v||\nabla q|+|q|^2|\nabla q|+|\nabla q|\right)|\nabla\d_t\,d\x\\
\leq&\int |\nabla\d_t|^2\,d\x+C\left[\|q\|_{\infty}^2\|\Delta v\|_2^2+\|\nabla v\|_4^4+\|\nabla q\|_4^4+(\|q\|_{\infty}^4+1)\|\nabla q\|_{2}^2\right]\\
\leq &\int |\nabla\d_t|^2\,d\x+C\left[(\| q\|_{H^2}^2+\|\nabla v\|_2^2)\|\Delta v\|_2^2+\|\nabla q\|_2^2\|\Delta q\|_2^2+(\|q\|_{\infty}^4+1)\|\nabla q\|_{2}^2\right]\\
\leq &\int |\nabla\d_t|^2\,d\x+CK^4+C(K^6+K^2)+C\left(K^2+\|\nabla v\|_2^2\right)\|\Delta v\|_2^2.
\end{split}
\eeq
Combining \eqref{prop1pf4}-\eqref{prop1pf5}, integrating over $(0,t)$ with respect to time and using the estimate \eqref{prop1pf31}, we have
\beq\label{prop1pf6}
\begin{split}
&\sup\limits_{0\leq s\leq t} \int \left(|\nabla \d_t|^2+|\Delta \d|^2\right)\,d\x +\int_0^t\int |\nabla \d_t|^2\,d\x ds+\delta\int_0^t\int |\Delta \d_t|^2\,d\x ds\\
\leq &\|\nabla \d_1\|_2^2+\|\Delta \d_0\|_2^2+Ct+C\int_0^t\|\Delta v\|_2^2\,ds.
\end{split}
\eeq
Similarly,  multiplying both sides of the second equation of \eqref{fels4}  by $ \d_t$, and integrating by parts over $\R^2\times (0,t)$, we have 
\beq\label{prop1pf41}
\begin{split}
\int \left(| \d_t|^2+|\nabla \d|^2\right)\,d\x +\int_0^t\int |\d_t|^2\,d\x ds+\delta\int_0^t\int |\nabla\d_t|^2\,d\x ds
\leq \| \d_0\|_2^2+\| \d_1\|_2^2+Ct.
\end{split}
\eeq
By this estimate, we have 
\beq\label{prop1pf42}
\begin{split}
\int |\d-\e_1|^2\,d\x\leq C\int |\nabla \d|^2\,d\x 
\leq C\| \d_0\|_2^2+C\| \d_1\|_2^2+Ct.
\end{split}
\eeq

Combining it with \eqref{prop1pf6}, \eqref{prop1pf41} and \eqref{prop1pf42}, we conclude 
\beq\notag
\|\d-\e_1\|_{\mathcal X}\leq C(\d_0, \d_1)+C(K)t+C\int_0^t\|\Delta v\|_2^2\,ds.
\eeq
If we take $K=2C(\d_0, \d_1)$ and choose $t$ small enough,  then $\d-\e_1\in \mathcal X_K$, which completes the proof of the claim. 

\medskip
\noindent{\it Claim 3. $\mathcal L$ is a contraction map with respect to distance $\rho_{\mathcal X}$ .} 

\medskip
\noindent Indeed, for any $q_i\in \mathcal X_K$ for $i=1,2$, let $v_i$ be solutions corresponding to $q_i$ and 
$$
\d_i=\mathcal L(q_i).
$$ 
Denote 
$$
v=v_1-v_2,\quad \d=\d_1-\d_2.
$$
Then $v$ and $\d$ satisfy the following equations with zero initial conditions 
\begin{equation}\label{fels5}
\begin{cases}
 v_t-\frac12\Delta v=\nabla\cdot\mathcal T_1(q_1,(q_1)_t)-\nabla\cdot\mathcal T_1(q_2,(q_2)_t)+\nabla\cdot\mathcal T_2(q_1,\nabla v_1)-\nabla\cdot\mathcal T_2(q_2,\nabla v_2),
\\
\d_{tt}-\delta\Delta\d_t+2 \d_t-\Delta \d=2\omega(v_1) q_1-2\omega(v_2)q_2-(|q_1|^2-1)q_1+(|q_2|^2-1)q_2.
\end{cases}
\end{equation} 
Multiplying the first equation of \eqref{fels5} by $v$ and integrating by parts over $\R^2$, we have 
\beq\label{prop1pf7}
\begin{split}
&\frac{d}{dt} \frac12\int |v|^2\,d\x +\frac12\int |\nabla v|^2\,d\x\\
=& -\int \left(T_1(q_1,(q_1)_t)-\mathcal T_1(q_2,(q_2)_t)\right)\cdot\nabla v\,d\x
-\int\left(\mathcal T_2(q_1,\nabla v_1)-\mathcal T_2(q_2,\nabla v_2)\right)\cdot\nabla v \,d\x\\
=&-\int \left(T_1(q_1,(q_1)_t)-\mathcal T_1(q_2,(q_2)_t)\right)\cdot\nabla v\,d\x
-\int\left(T_2(q_1,\nabla v_2)-\mathcal T_2(q_2,\nabla v_2)\right)\cdot\nabla v  \,d\x\\
&-\int\mathcal T_2(q_1,\nabla v)\cdot\nabla v\,d\x.
\end{split}
\eeq
For the first term, it holds
\beq\label{prop1pf8}
\begin{split}
&\left|\int \left(T_1(q_1,(q_1)_t)-\mathcal T_1(q_2,(q_2)_t)\right)\cdot\nabla v\,d\x\right|\\
\leq&C\int\left(|q||(q_2)_t||\nabla v|+|q_2||q_t||\nabla v|\right)\,d\x\\
\leq &\frac18\int|\nabla v|^2\,d\x+C\left(\|q\|_4^2\|(q_2)_t\|_4^2+\|q_2\|_\infty^2\|q_t\|_2^2\right)\\
\leq &\frac18\int|\nabla v|^2\,d\x+CK\left(\|q\|_2\|\nabla q\|_2+\|q_t\|_2^2\right)\\
\leq &\frac18\int|\nabla v|^2\,d\x+C(K)\left(\|q\|_2^2+\|\nabla q\|^2_2+\|q_t\|_2^2\right).
\end{split}
\eeq
Similarly, it holds for the second term
 \beq\label{prop1pf9}
\begin{split}
&\left|\int\left(T_2(q_1,\nabla v_2)-\mathcal T_2(q_2,\nabla v_2)\right) \,d\x\right|\\
\leq&C\int|q||q_1+q_2||\nabla v_2||\nabla v|\,d\x\\
\leq &\frac18\int|\nabla v|^2\,d\x+C\left(\|q_1\|_\infty^2+\|q_2\|_\infty^2\right)\|q\|_4^2\|\nabla v_2\|_4^2\\
\leq &\frac18\int|\nabla v|^2\,d\x+C(K)\|q\|_2\|\nabla q\|_2\|\nabla v_2\|_2\|\Delta v_2\|_2\\
\leq &\frac18\int|\nabla v|^2\,d\x+C(K)\|\Delta v_2\|_2(\|q\|_2^2+\|\nabla q\|_2^2).
\end{split}
\eeq
By the computation in \eqref{lem1pf3}, we have
\beq\label{prop1pf10}
\begin{split}
-\int\mathcal T_2(q_1,\nabla p)\cdot\nabla v\,d\x\leq 0.
\end{split}
\eeq
Combining \eqref{prop1pf7}-\eqref{prop1pf10} and integrating over $(0,t)$, it holds
\beq\label{prop1pf11}
\begin{split}
&\sup\limits_{0\leq s\leq t}\int |v|^2\,d\x +\int_0^t\int |\nabla v|^2\,d\x ds\\
\leq& \left(Ct+C\int_0^t\|\Delta v\|_2\,ds\right)\rho^2_{\mathcal X}(q_1,q_2)\\
\leq &C(K)(t+\sqrt{t})\rho^2_{\mathcal X}(q_1,q_2).
\end{split}
\eeq
Multiplying the second equation of \eqref{fels5} by $\d_t$ and integrating by parts over $\R^2$, we have 
\beq\notag
\begin{split}
&\frac{d}{dt} \frac12\int |\d_t|^2\,d\x +\frac{d}{dt} \frac12\int |\nabla \d|^2\,d\x+2\int |\d_t|^2\,d\x
+\delta\int |\nabla\d_t|^2\,d\x\\
=& \int\big(2\omega(v_1) q_1-2\omega(v_2)q_2-(|q_1|^2-1)q_1+(|q_2|^2-1)q_2\big)\d_t\,d\x\\
\leq& C\int\big(|\nabla v||q_1|+|\nabla v_2||q|+|q|(|q_1|^2+|q_2|^2)+|q|\big)|\d_t|\,d\x\\
\leq & \frac12\int |\d_t|^2\,d\x+C\left(\|q_1\|_{\infty}^2\|\nabla v\|_2^2+\|q\|_{4}^2\|\nabla v_2\|_4^2\right)+C\int|q|^2\big(|q_1|^4+|q_2|^4+1\big)\,d\x,
\end{split}
\eeq
which implies 
\beq\notag
\begin{split}
&\frac{d}{dt}\int |\d_t|^2\,d\x +\frac{d}{dt} \int |\nabla \d|^2\,d\x+\int |\d_t|^2\,d\x
+\delta\int |\nabla\d_t|^2\,d\x\\
\leq &C\left(\|q_1\|_{H^2}^2\|\nabla v\|_2^2+\|q\|_{2}\|\nabla q\|_2\|\nabla v_2\|_2\|\Delta v_2\|_2\right)+C\big(\|q_1\|_{\infty}^4+\|q_2\|_{\infty}^4+1\big)\|q\|_2^2\\
\leq &C\|\nabla v\|_2^2+C\|\Delta v_2\|_2(\|q\|_{2}^2+\|\nabla q\|_2^2)+C\|q\|_2^2.
\end{split}
\eeq
Integrating over $(0,t)$ and using \eqref{prop1pf11}, we conclude
\beq\label{prop1pf12}
\rho^2_{\mathcal X}(\d_1,\d_2)=\sup\limits_{0\leq s\leq t}\left(\int |\d_t|^2\,d\x + \int |\nabla \d|^2\,d\x\right)
\leq C(K)(t+\sqrt{t})\rho^2_{\mathcal X}(q_1,q_2).
\eeq
Choosing $t$ sufficient small, we prove the claim. The existence of the solutions to system \eqref{felsd} can now be derived by the standard fixed-point argument providing initial conditions \eqref{initial5}. 

\medskip
\noindent{\it Claim 4. The solutions obtained in Claim 3 is global.} 

\medskip
\noindent
Similar to the proof of Lemma \ref{lemma2.1}, we can prove the energy estimate for the system \eqref{felsd}.
\beq\label{engestwd}
\begin{split}
&\int \left(|v|^2+ | \d_t|^2+|\nabla \d|^2+\left(|\d|^2-1\right)^2\right)\,d\x
+C\int\int |\nabla v|^2\,d\x ds+\delta\int\int |\nabla\d_t|^2\,d\x ds\\
\leq& \int\left(|v_0|^2+ | \d_1|^2+|\nabla \d_0|^2\right)\,d\x
\end{split}
\eeq
It suffices to show the uniform estimate for $\|\nabla \d_t\|_2^2+\|\Delta \d\|_2^2$ for any $t>0$. Multiplying the second equation of \eqref{felsd} by $\Delta d_t$ and integrating over $\R^2$, we have 
\beq\notag
\begin{split}
&\frac{d}{dt} \frac12\int \left(|\nabla \d_t|^2+|\Delta \d|^2\right)\,d\x +2\int |\nabla\d_t|^2\,d\x+\delta\int|\Delta \d_t|^2\,d\x\\
=&-\int 2\omega(v) \d\cdot \Delta\d_t\,d\x+\int(|\d|^2-1)\d\cdot \Delta\d_t\,d\x\\
\leq &\frac{\delta}{2}\int|\Delta \d_t|^2\,d\x+C\int\left(|\nabla v|^2|\d|^2+(|\d|^2-1)^2|\d|^2\right)\,dx
\end{split}
\eeq
which implies
\beq\label{prop1pf51}
\begin{split}
&\frac{d}{dt} \frac12\int \left(|\nabla \d_t|^2+|\Delta \d|^2\right)\,d\x +2\int |\nabla\d_t|^2\,d\x+\frac{\delta}{2}\int|\Delta \d_t|^2\,d\x\\
\leq& C\int\left(|\nabla v|^2|\d|^2+(|\d|^2-1)^2|\d|^2\right)\,dx\\
\leq & C(\|\d-\e_1\|_{\infty}^2+1)\int\left(|\nabla v|^2+(|\d|^2-1)^2\right)\,dx\\
\leq & C\left(\int|\nabla v|^2\,dx+C\right)(\|\Delta \d\|_2^2+1).
\end{split}
\eeq
Standard Gronwall type estimate shows $\|\nabla \d_t\|_2^2+\|\Delta \d\|_2^2$ is bounded for any $0<t<\infty$, which completes the proof of the claim.

\medskip
For general initial conditions \eqref{initial4}, we can show the existence of weak solutions to the system \eqref{fels3} by letting $\delta\rightarrow 0^+$ and standard approximations of initial data $\d_0^l\in H^2$, $\d_1^l\in H^1$ for $l=1,2...$ such that 
\beq\notag
\|\d_0^l-\d_0\|_{H^1}+\|\d_1^l-\d_1\|_{L^2}\rightarrow 0
\eeq
as $l\rightarrow \infty$. By the uniform energy estimates \eqref{engestwd}, we may only get weak solutions. 
\endpf

\bigskip
As a result of Proposition \ref{prop2.3}, we can take the limit in \eqref{fels3} by the similar arguments as in \cite{Shatah01}.

Since $\d^{\varepsilon}$ is uniformly bounded in $L^\infty_tH^1_x$, and satisfies the energy equality \eqref{engestw2} for each $\varepsilon$, let $\varepsilon\to 0$, then there is a subsequence, still denoted as $(v^\varepsilon, \d^{\varepsilon})$, such that
\[
\d^\varepsilon\to \d, ~~W^* \text{ in } L^\infty_{loc}(\R, \dot H^1),\quad \partial_t\d^\varepsilon \to \d^{(t)},~~ W^* \text{ in } L^\infty_{loc}(\R, L^2)
\]
and 
\[
\int (|\d|^2-1)^2 d\x\leq\lim\limits_{\varepsilon\to0}\int (|\d^\varepsilon|^2-1)^2 d\x=0 \Rightarrow |\d|=1.
\]
By compactness and a diagonalization process we can choose the subsequence
$\{ \d^\varepsilon\}$ to converge strongly in $L^2_{loc}$ and a.e.. Moreover, by multiplying $\partial_t\d^{\varepsilon}$ by test functions we obtain $\partial_t\d=\d^{(t)}$.

In order to show that $\d$ satisfies the equation we take the cross product of  \eqref{fels3}$_2$ with $\d^\varepsilon$, to obtain
\[
\partial_t(\d^\varepsilon_{t}\times\d^\varepsilon)+2\d^\varepsilon_t\times\d^\varepsilon-2\omega\d^\varepsilon\times\d^\varepsilon=\nabla\cdot(\nabla\d^\varepsilon\times\d^\varepsilon).
\]
It satisfies 
\[
\iint -\partial_t\varphi(\d^\varepsilon_{t}\times\d^\varepsilon)+2\varphi\d^\varepsilon_t\times\d^\varepsilon-2\varphi\omega\d^\varepsilon\times\d^\varepsilon+\nabla\varphi\cdot(\nabla\d^\varepsilon\times\d^\varepsilon)d\x dt=0.
\]
Since $\d^{\varepsilon}$ converges strongly in $L^2_{loc}$ and $\partial \d^{\varepsilon}$ weakly in $L^2_{loc}$, we can pass to the limit to obtain
\[
\partial_t(\d_{t}\times\d)+2\d_t\times\d-2\omega\d\times\d=\nabla\cdot(\nabla\d\times\d),
\]
which implies
\[
\d_{tt}\times\d+2\d_t\times\d-2\omega\d\times\d=\Delta\d\times\d.
\]
Considering the fact that $|\d|=1$, and $\omega$ is anti-symmetric, we have
\[
\partial_t(\d_t\cdot\d)+2\d_t\cdot\d-2(\omega\d)\cdot \d=\nabla\cdot(\nabla\d\cdot\d).
\]
Therefore, combining the above two equations, we obtain \eqref{fels2}$_2$.

Since $v^\varepsilon$ satisfies \eqref{fels3}$_1$,
\[
\iint -v^\varepsilon\varphi_t+\frac12\nabla v^\varepsilon\cdot\nabla\varphi+\mathcal T_1(\d^\varepsilon, \d^\varepsilon_t)\cdot\nabla\varphi+\mathcal T_2(\d^\varepsilon, \nabla v^\varepsilon)\cdot\nabla\varphi d\x dt=0.
\]
Notice that $v^\varepsilon$ and $\nabla v^\varepsilon$ converge weakly in $L^2((0,T),L^2(\R^2))$, therefore
\[
\iint -v^\varepsilon\varphi_t+\frac12\nabla v^\varepsilon\cdot\nabla\varphi d\x dt\to \iint -v\varphi_t+\frac12\nabla v\cdot\nabla\varphi d\x dt.
\]
For $\mathcal T_1$ and $\mathcal T_2$, we have
\[
\mathcal T_1(\d^\varepsilon, \d^\varepsilon_t)+\mathcal T_2(\d^\varepsilon, \nabla v^\varepsilon)=
\left(\begin{array}{c}
d^{3,\varepsilon}(d^{1,\varepsilon}_{t}+\frac12 v_xd^{3,\varepsilon})-d^{1,\varepsilon}(d^{3,\varepsilon}_{t}-\frac12v_xd^{1,\varepsilon}-\frac12 v_yd^{2,\varepsilon})\\
d^{3,\varepsilon}(d^{2,\varepsilon}_{t}+\frac12 v_xd^{3,\varepsilon})-d^{2,\varepsilon}(d^{3,\varepsilon}_{t}-\frac12v_xd^{2,\varepsilon}-\frac12 v_yd^{1,\varepsilon})
\end{array}
\right).
\]
By the energy estimate Lemma \ref{lemma2.2}, we have $\d^\varepsilon_t-\omega\d^\varepsilon\in L^2((0,T), L^2(\R^2))$ 
is bounded uniformly with respect to $\varepsilon$, thus weakly convergent. Since $\d^{\varepsilon}$ converges strongly in $L^2_{loc}$, $\mathcal T_1+\mathcal T_2$ weakly converges in $L^2_{loc}$. 

Since $\d^\varepsilon_t-\omega\d^\varepsilon$ converges weakly in $L^2((0,T),L^2_{loc})$, and  $\d^{\varepsilon}$ converges strongly in $L^2_{loc}$,
we obtain 
\[
\iint\mathcal T_1(\d^\varepsilon, \d^\varepsilon_t)\cdot\nabla\varphi+\mathcal T_2(\d^\varepsilon, \nabla v^\varepsilon)\cdot\nabla\varphi d\x dt\to  \iint\mathcal T_1(\d, \d_t)\cdot\nabla\varphi+\mathcal T_2(\d, \nabla v)\cdot\nabla\varphi d\x dt.
\]
Therefore we obtain the existence of weak solutions to \eqref{fels2}.

In addition, if the initial data is axisymmetric, then the approximate solutions $(v^\varepsilon, \d^\varepsilon)$ are axisymmetric, and there limit is also axisymmetric. So we proved the following proposition.

\begin{proposition}\label{prop2.4}
For any $t>0$, the Cauchy problem of \eqref{fels2} with initial and boundary conditions \eqref{initial2}-\eqref{initial4} has a weak solution 
\beq\notag
v\in L^\infty((0,t), L^2(\R^2,\R))\cap L^2((0,t), H^1(\R^2,\R)),
\eeq
\beq\notag
\d\in L^\infty((0,t), H^1(\R^2,\mathbb S^2)),\quad  \d_t\in L^\infty((0,t), L^2(\R^2,\R^3)).
\eeq 
In additionally, the solutions $(v,\d)$ depend only on $(r,t)$ and $\d$ is axisymmetric without swirl, which implies the existence of global weak solutions to the system \eqref{wels} with initial condition \eqref{welsr}.
\end{proposition}


\section{Energy estimates}
\setcounter{equation}{0}

This section will be devoted to energy estimates which are useful in proof of partial regularity. 
The first estimate is for the system \eqref{wels} obtained in \cite{CHX23}
\begin{lemma}\label{englemma1}
For any smooth solution to the system \eqref{wels}, we have the energy estimate,
\beq\label{enginq1}
\begin{split}
\frac{d}{dt}\frac{1}{2}\int(|v|^2+|\phi_r|^2+|\phi_t|^2+\frac{\sin^2(\phi)}{r^2})rdr+c\int(|\phi_t|^2+|v_r|^2)rdr\leq 0,
\end{split}
\eeq
for some constant $c>0$.
\end{lemma}

To obtain better energy estimates, we need to utilize the system \eqref{welss} of $(h,\phi)$. 

\begin{lemma}\label{englemma2}
For any smooth solution to \eqref{welss}, with boundary conditions \eqref{welsr} and \eqref{welshr}, we have
\beq\label{enginq22}
\begin{split}
\frac{d}{dt}\frac{1}{2}\int\left(|\phi_r|^2+|\phi_t|^2+\frac{\sin^2(\phi)}{r^2}+|h_r|^2+\frac{h^2}{r^2}\right)\,rdr+\int\left(|\phi_t|^2+|h_t|^2\right)\,rdr\leq 0.
\end{split}
\eeq
\end{lemma}

\pf Multiplying the first equation of \eqref{welss} by $h_tr$ and the second equation by $\phi_tr$, adding two resulting equations and integrating over $(0,\infty)$ with respect to $r$, we obtain the energy estimate \eqref{enginq22}.
\endpf

We can also derive the higher energy estimate for \eqref{welss} as in \cite{CHX23}.

\begin{lemma}\label{englemma3}
For any smooth solution to \eqref{welss}, with boundary conditions \eqref{welsr} and \eqref{welshr}, we have
\beq\label{enginq2}
\begin{split}
\frac{d}{dt}\frac{1}{2}\int|h_t|^2\, rdr+\int\left(|h_t|^2+\frac{|h_t|^2}{r^2}+|h_{tr}|^2\right)\,rdr\leq CE_0.
\end{split}
\eeq
Furthermore, if $E_0<\infty$, it holds
\beq\label{enginq3}
\begin{split}
\sup_{t\ge0}\int|h_t|^2rdr\leq \int|h_t|^2(r,0)rdr+Ct.
\end{split}
\eeq
\end{lemma}

\begin{remark}\label{remark1}
By the definition of $h(r,t)$, it is not hard to see that 
$$
h_t=\left(\frac{1}{r}\int_0^rv(R,t)\,RdR\right)_t=v_r+\phi_t=v_r+\phi_t
$$
which implies that $h_t(r,0)$ is well-defined in $L^2$.
\end{remark}

\medskip

\begin{lemma}\label{englemma4}
We have the following local energy estimate for any $0< \tau \leq s< T$,
\beq\label{eneflux1}
\begin{split}
E(R,s)+\text{Flux}(s,s-\tau)\leq E(R+\tau,s-\tau)+C\tau.
\end{split}
\eeq
In particular, we also have the following estimates for any $R+s=T$
\beq\label{eneflux2}
\begin{split}
E(R,s)\leq E(T,0)+CT.
\end{split}
\eeq
\end{lemma}
\pf
Multiplying the second equation of \eqref{welss} by $\phi_tr$ we obtain
\beq\label{eng3}
\begin{split}
\partial_t\left(\frac{1}{2}\phi_t^2r+\frac{1}{2}\phi_r^2r+\frac{1}{2}\frac{\sin^2\phi}{r^2}r\right)+\phi^2_tr-\partial_r(\phi_r\phi_tr)=-h_t\phi_tr.
\end{split}
\eeq
Integrating \eqref{eng3} over $K^T(s,\tau)$, we have
\beq\notag
\begin{split}
\iint_{K^T(s,\tau)}\big((er)_t-(\phi_r\phi_tr)_r\big)\,drdt+\iint_{K^T(s,\tau)}|\phi_t|^2\,rdrdt=-\iint_{K^T(s,\tau)}h_t\phi_t\,rdrdt.
\end{split}
\eeq
For the 1st term, by Green's theorem, it holds
\beq\notag
\begin{split}
&\iint_{K^T(s,\tau)}[(er)_t-(\phi_r\phi_tr)_r]\,drdt\\
=&-\int_{M^T(s,\tau)}\phi_r\phi_t\,rdt+e\,rdr+E(R,s)-E(R+\tau,s-\tau)\\
=&\text{Flux}(s,s-\tau)+E(R,s)-E(R+\tau,s-\tau)
\end{split}
\eeq
Combining these estimates together and using Young's inequality, we have
\beq\notag
\begin{split}
&E(R,s)-+\text{Flux}(s,s-\tau)+\frac{1}{2}\iint_{K^T(s,\tau)}|\phi_t|^2\,rdrdt\\
\leq& E(R+\tau,s-\tau)+ C\int_{s-\tau}^s\int_0^{T-t}|h_t|^2\,rdrdt\\
\leq &E(R+\tau,s-\tau)+C\tau, 
\end{split}
\eeq
where we have used the Remark \ref{remark1} and the energy estimate \eqref{enginq3} in last inequality.
This completes the proof of lemma.
\endpf
\begin{lemma}\label{englemma5}
It holds that
\beq\label{flux0}
\begin{split}
\text{Flux}\,(T,T-\tau)\rightarrow 0 \text{ as } \tau\rightarrow0^+. 
\end{split}
\eeq
\end{lemma}
\pf
Combining the estimates \eqref{eneflux1} and \eqref{eneflux1} with the fact $R+s=T$, we obtain 
\beq\notag
\begin{split}
\text{Flux}\,(s,s-\tau)
\leq E(R+\tau,s-\tau)+C\tau
\leq E(T,0)+CT
\leq E_0+CT. 
\end{split}
\eeq
This implies that for any $0\leq \tau$
$$ \text{Flux}\,(T,T-\tau)\leq E_0+CT,$$
and furthermore
$$\text{Flux}\,(T,T-\tau)=\int_{T-\tau}^T\frac{1}{2}(\chi')^2+\frac{F(\chi)}{l^2}dl\rightarrow 0 \text{ as } \tau\rightarrow 0^+$$
\endpf
\begin{lemma}\label{englemma6}
For any $0<t<T$ and $R>0$, we have
\beq
\begin{split}\label{upb}
\sup_{r<R}|\phi(r,t)|\leq C(E(R,t))
\end{split}
\eeq
\end{lemma}
\pf
Denote $$H(\phi(r,t))=\int_0^\phi|\sin \tau|\,d\tau.$$
The fact $H(0)=0$ implies
\beq\notag
\begin{split}
H(\phi)=\int_0^r\frac{\partial}{\partial r}H(\phi(r,t))\,dr=\int_0^r|\sin\phi|\phi_r\,dr.
\end{split}
\eeq
Hence
\beq\notag
\begin{split}
 |H(\phi)|\leq \left(\int_0^r\frac{\sin^2\phi}{r^2}\,rdr\right)^\frac{1}{2}\left(\int_0^r|\phi_r|^2\,rdr\right)^\frac{1}{2}\leq CE(R,t).
\end{split}
\eeq
Since $|H(\phi)|\rightarrow \infty$ as $|\phi|\rightarrow\infty$, we conclude \eqref{upb}, which completes the proof of lemma.
\endpf

\begin{lemma}\label{englemma7} For any $0<\lambda<1$, it holds
\beq\label{ext0}
\lim_{t\rightarrow T^-}\int_{\lambda(T-t)}^{T-t}e(r,t)\,rdr
=\lim_{\tau\rightarrow 0^+}\int_{\lambda\tau}^\tau e(r,T-t)\,rdr
=0.
\eeq
\end{lemma}

\pf
Denote $m=\phi_t\phi_r$ and recall that $e=\frac{1}{2}\phi_t^2+\frac{1}{2}\phi_r^2+\frac{1}{2}\frac{\sin^2\phi}{r^2}$.
By the equation \eqref{eng3} and multiplying the second equation  of \eqref{welss} by $r\phi_r$, we obtain
\begin{equation}\label{wels3}
\begin{cases}
(re)_t-(rm)_r+\phi_t^2r=-h_t\phi_tr\\
(rm)_t-(re)_r+r\phi_r\phi_t=L-h_t\phi_rr,
\end{cases}
\end{equation}
where
$$
L=-\frac{1}{2}\phi_t^2+\frac{1}{2}\phi_r^2+\frac{1}{2}\frac{\sin^2\phi}{r^2}-\frac{\sin(2\phi)}{r}\phi_r.
$$
Let's introduce the coordinate changes
\begin{equation}\label{tran1}
\begin{cases}
\xi=t+r\\
\eta=t-r
\end{cases}
\quad 
\begin{cases}
t=\frac{\xi+\eta}{2}\\
r=\frac{\xi-\eta}{2}.
\end{cases}
\end{equation}
Then it is easy to see  
$$\partial_\xi=\frac{1}{2}\partial_t+\frac{1}{2}\partial_r,\partial_\eta=\frac{1}{2}\partial_t-\frac{1}{2}\partial_r.$$
Denote
$$
\alpha^2=r(e-m)\ge0,\quad 
\beta^2=r(e+m)\ge0
$$
The system \eqref{wels3} becomes
$$
\partial_\xi\alpha^2=\frac{1}{2}[-\phi_t^2r-h_t\phi_tr-L+\phi_r\phi_tr+h_t\phi_rr]
$$
$$
\partial_\eta\beta^2=\frac{1}{2}[-\phi_t^2r-h_t\phi_tr+L-\phi_r\phi_tr-h_t\phi_rr].
$$
Direct computation implies
$$
L^2\leq\frac{c}{r^2}\alpha^2\beta^2
$$
$$
\frac{\alpha^2}{r}=e-m=\frac{1}{2}(\phi_t-\phi_r)^2+\frac{1}{2}\frac{\sin^2\phi}{r^2}
$$
$$
\frac{\beta^2}{r}=e+m=\frac{1}{2}(\phi_t+\phi_r)^2+\frac{1}{2}\frac{\sin^2\phi}{r^2}
$$
$$
\frac{1}{2}[-\phi_t^2r-h_t\phi_tr+\phi_r\phi_tr+h_t\phi_rr]\leq \frac{\sqrt{r}}{2}\alpha(|\phi_t|+|h_t|)
$$
$$
\frac{1}{2}[-\phi_t^2r-h_t\phi_tr-\phi_r\phi_tr-h_t\phi_rr]\leq \frac{\sqrt{r}}{2}\beta(|\phi_t|+|h_t|)
$$
Combining these estimates with the equations of $\alpha$ and $\beta$, we obtain
\begin{equation}\label{wel4}
|\partial_\xi\alpha|\leq \frac{c}{r}\beta+c\sqrt{r}(|\phi_t|+|h_t|),\quad 
|\partial_\eta\beta|\leq \frac{c}{r}\alpha+c\sqrt{r}(|\phi_t|+|h_t|).
\end{equation}
Let's consider a region $S(\xi,\eta)=\{(\xi',\eta')\,\big|\,\xi\leq\xi'\leq T,\quad\eta_0\leq\eta'\leq\eta\}$ as a subset of 
$
K^T_{\lambda}= \{0\leq t\leq T,\ \lambda(T-t)\leq r\leq T-t\}.
$
By the definition of $\xi$ and $\eta$, it is easy to see $\xi\ge\eta\ge\eta_0$. 
Direct computation implies 
\beq\label{TT}
\begin{split}
&\beta(\xi,\eta)-\beta(\xi,\eta_0)\\
\lesssim&\int_{\eta_0}^\eta\frac{\alpha}{\xi-\eta'}(T,\eta')d\eta'+\int_{\eta_0}^\eta\sqrt{r}(|\phi_t|+|h_t|)d\eta'+\int_{\eta_0}^\eta\int_\xi^{T} \frac{\beta(\xi',\eta')}{(\xi-\eta')(\xi'-\eta')}d\xi'd\eta'\\
&+\int_{\eta_0}^\eta\int_\xi^{T}\sqrt{r}(|\phi_t|+|h_t|)d\xi'd\eta' \triangleq L_1+L_2+L_3+L_4.
\end{split}
\eeq
The first term $L_1$ can be estimated as follows
\beq\notag
\begin{split}
|L_1|&\leq(\int^\eta_{\eta_0}\alpha^2(T,\eta')d\eta')^\frac{1}{2}(\int^\eta_{\eta_0}\frac{d\eta'}{(\xi-\eta')^2})^\frac{1}{2}\\
&=(\int^\eta_{\eta_0}r(e-m)(T,\eta')d\eta')^\frac{1}{2}(\frac{1}{\xi-\eta}-\frac{1}{\xi-\eta_0})^\frac{1}{2} \\
&\lesssim(\text{Flux}\,(T,\eta_0))^\frac{1}{2}/(\xi-\eta)^\frac{1}{2}.
\end{split}
\eeq
The terms $L_2,\ L_4$ can be estimated as follows
\beq\notag
\begin{split}
|L_2|\leq\left(\int^\eta_{\eta_0}(|\phi_t|^2+|h_t|^2)\,rd\eta'\right)^\frac{1}{2}\left(\eta-\eta_0\right)^\frac{1}{2},
\end{split}
\eeq
\beq\notag
\begin{split}
|L_4|\leq\left(\int^\eta_{\eta_0}\int^T_\xi(|\phi_t|^2+|h_t|^2)\,rd\xi'd\eta'\right)^\frac{1}{2}\text{Area}( S(\xi,\eta))^\frac{1}{2}
\lesssim \text{Area}(S(\xi,\eta))^\frac{1}{2},
\end{split}
\eeq
where we have the energy estimate \eqref{enginq22} in last inequality.
For $L_3$:
\beq\notag
\begin{split}
|L_3|&\leq \left(\int^\eta_{\eta_0}\int^T_\xi \frac{\beta^2(\xi',\eta')}{(\xi-\eta')^2}\,d\xi'd\eta'\right)^\frac{1}{2}\left(\int^\eta_{\eta_0}\int^T_\xi\frac{1}{(\xi'-\eta')^2}\,d\xi'd\eta'\right)^\frac{1}{2}\\
&=\left(\int^\eta_{\eta_0}\int^T_\xi\frac{\beta^2(\xi',\eta')}{(\xi-\eta')^2}\,d\xi'd\eta'\right)^\frac{1}{2}\left(\int^T_\xi(\frac{1}{\xi'-\eta}-\frac{1}{\xi'-\eta_0})\,d\xi'\right)^\frac{1}{2}\\
&\leq\left(\int^\eta_{\eta_0}\int^T_\xi\frac{\beta^2(\xi',\eta')}{(\xi-\eta')^2}\,d\xi'd\eta'\right)^\frac{1}{2}\left(\ln{\frac{T-\eta}{\xi-\eta}}\right)^\frac{1}{2}.
\end{split}
\eeq
In case of $(\xi,\eta)$ is on the line $r+\lambda t=\lambda T$, it holds
$$
\xi-\eta=\frac{2\lambda}{1+\lambda}(T-\eta)\Rightarrow \ln{\frac{T-\eta}{\xi-\eta}}=\ln{\frac{2\lambda}{1+\lambda}}.
$$
Hence, we may choose $(\xi,\eta)$ close to $(T,T)$ such that $\ln{\frac{T-\eta}{\xi-\eta}}\leq\ln{(\frac{2\lambda}{1+\lambda})}$, which implies
$$
 |L_3|\lesssim\left(\int^\eta_{\eta_0}\int^T_\xi\frac{\eta^2(\xi',\eta')}{(\xi-\eta')^2}\,d\xi'd\eta'\right)^\frac{1}{2}.
$$
Plugging all these estimates into \eqref{TT}, squaring both sides of resulting inequality and integrating over $(\xi,T)$, we derive
\beq\label{betaest1}
\begin{split}
\int^T_\xi \beta^2(\xi',\eta)\,d\xi'\lesssim &\int^T_\xi\beta^2(\xi',\eta_0)\,d\xi'+
\int^T_\xi\frac{\text{Flux}\,(T,\eta_0)}{\xi'-\eta}\,d\xi'\\
&+(\eta-\eta_0)\int^T_\xi\int^\eta_{\eta_0}(|\phi_t|^2+|h_t|^2)rd\eta'\,d\xi'+\text{Area}(S(\xi,\eta))(T-\xi)\\&+\int^T_\xi\int^\eta_{\eta_0}\int^T_{\xi'}\frac{\beta^2(\xi'',\eta')}{(\xi'-\eta')^2}\,d\xi''d\eta'd\xi'\\&\triangleq M_1+M_2+M_3+M_4
\end{split}
\eeq
Since the solution is regular before time $T$, we have
$$
|M_1|\leq C(\eta_0)(T-\xi).
$$
The rest terms can be estimated as follows
$$
|M_2|=\text{Flux}(T,\eta_0)\ln{\frac{T-\eta}{\xi-\eta}}\leq\text{Flux}(T,\eta_0)\rightarrow 0 \text{ as }\eta_0\rightarrow T^-,
$$
$$
|M_3|\lesssim C(T-\xi)\rightarrow 0 \text{ as } \xi\rightarrow T^-,
$$
\beq\notag 
\begin{split}
M_4&=\int^\eta_{\eta_0}\int^T_\xi\int^T_{\xi'}\frac{\beta^2(\xi'',\eta')}{(\xi'-\eta')^2}\,d\xi''d\xi'd\eta'
\\&=\int^\eta_{\eta_0}\int^T_\xi\int^{\xi''}_\xi\frac{\beta^2(\eta'',\xi')}{(\xi'-\eta')^2}\,d\xi'd\xi''d\eta'
\\&=\int^\eta_{\eta_0}\int^T_\xi\beta^2(\xi'',\eta')\left(\frac{1}{\xi-\eta'}-\frac{1}{\xi''-\eta'}\right)\,d\xi''d\eta'
\\&\leq\int^\eta_{\eta_0}\left(\frac{1}{\eta-\xi'}-\frac{1}{T-\eta'}\right)\int^T_\xi\beta^2(\xi'',\eta')\,d\xi''d\eta'.
\end{split}
\eeq
Denote 
$$
F(\xi,\eta)=\int_\xi^T\beta^2(\xi',\eta)\,d\xi',\quad M_0=C(\eta_0)(T-\xi)+C\text{Flux}(T,\eta_0)
$$
Combining all the above estimates with \eqref{betaest1}, we obtain 
$$
F(\xi,\eta)\leq M_0+\int^\eta_{\eta_0}\left(\frac{1}{\xi-\eta'}-\frac{1}{T-\xi}\right)F(\xi,\eta')\,d\eta'.
$$
Gronwall's inequality implies
$$
F(\xi,\eta)\leq M_0+\int^\eta_{\eta_0}M_0\left(\frac{1}{\xi-\eta'}-\frac{1}{T-\eta'}\right)e^{\int_{\eta'}^\eta\left(\frac{1}{\xi-\eta''}-\frac{1}{\T-\eta''}\right)\,d\eta''}\,d\eta'.
$$
Direct computation implies
\beq\notag
\begin{split}
\int^\eta_{\eta'}\left(\frac{1}{\xi-\eta''}-\frac{1}{T-\eta''}\right)\,d\xi''
=\ln\frac{T-\eta}{\xi-\eta}-\ln\frac{T-\eta'}{\xi-\eta'}
\leq \ln\frac{2\lambda}{1+\lambda}.
\end{split}
\eeq
Therefore, we conclude 
\beq\label{best0}
 \int^T_\xi\beta^2(\xi',\eta)d\xi'=F(\xi,\eta)\lesssim M_0\rightarrow 0
 \eeq
 as $\xi,\ \eta_0,\ \eta\rightarrow T^-$. 

Now we are ready to prove the estimate \eqref{ext0}.
We consider the triangle region $\triangle$ with vertices $(\xi, \eta)$, $(T, \eta_1)$ and $(T, \eta)$ where $(\xi, \eta)$ is on the line segment $r+\lambda t=\lambda T$ with $0\leq t\leq T$ and $\eta_1=\xi+\eta-T$.
Integrating the first equation of \eqref{wels3} over the triangle region, we obtain
\beq\label{pf1}
\begin{split}
 \iint_\triangle(\phi^2_tr+h_t\phi_tr)\,drdt
=-\iint_\triangle (re)_t-(rm)_r)\,drdt
=\int_{S_1}e\,rdr+(\int_{S_2}+\int_{S_3})rm\,dt+re\,dr,
\end{split}
\eeq
where $S_i$ for $i=1,2,3$ are three sides of the triangle $\triangle$ starting from $(\xi,\eta)$ with counterclockwise orientation.
Direct computation implies
$$
\int_{S_2}rm\,dt+re\,dr=  \text{Flux}\,(T,(\xi+\eta)/2),
$$
and
$$
\int_{S_3}rm\,dt+re\,dr=\frac{1}{2}\int^T_\xi(rm+re)(\xi',\eta)\,d\xi'=\frac{1}{2}\int^T_\xi\beta^2(\xi',\eta)\,d\xi'.
$$
Therefore
\beq\notag
\int_{\lambda(T-t)}^{T-t}e(r,t)\,rdr=\int_{S_1}e\,rdr
= \iint_\triangle(\phi^2_tr+h_t\phi_tr)\,drdt-\text{Flux}\,(T,(\xi+\eta)/2)
-\frac{1}{2}\int^T_\xi\beta^2(\xi',\eta)\,d\xi'.
\eeq
The right side approaches to zero as $\xi, \eta\rightarrow T^-$ by estimates \eqref{enginq22}, \eqref{flux0} and \eqref{best0},
which completes the proof of lemma.
\endpf

\begin{lemma}\label{englemma8}
Denote $K^T(\tau)= \{T-\tau\leq t\leq T,\ 0\leq r\leq T-t\}$.
Then it holds
\beq\label{phit0}
\lim\limits_{\tau\rightarrow 0^+}\frac{1}{\tau}\iint_{K^T(\tau)}|\phi_t|^2\, rdrdt\rightarrow 0. 
\eeq
\end{lemma}

\pf
Multiplying both sides of the second equation of \eqref{welss} by $r^2\phi_r$, we obtain
\begin{equation}\notag
(r^2\phi_t\phi_r)_t-\left(\frac{1}{2}r^2\phi_t^2+\frac{1}{2}r^2\phi_r^2-\frac{1}{2}\sin^2\phi\right)_r+r^2\phi_t\phi_r+r^2h_th_r+r\phi_t^2=0.
\end{equation}
Integrating it over $K^T(s,\tau)$
\beq\label{pf2}
\begin{split}
\iint_{K^T(s,\tau)}\phi^2_t\,rdrdt=&-\iint_{K^T(s,\tau)}\left(r^2\phi_t\phi_r+r^2h_t\phi_r\right)\,drdt\\
&-\iint_{K^T(s,\tau)}\left((r^2\phi_t\phi_r)_t-\left(\frac{1}{2}r^2\phi_t^2+\frac{1}{2}r^2\phi_r^2-\frac{1}{2}\sin^2\phi\right)_r\right)\,drdt.
\end{split}
\eeq
For the fist term of right side, it holds
\beq\notag
\begin{split}
&\left|\iint_{K^T(s,\tau)}\left(r^2\phi_t\phi_r+r^2h_t\phi_r\right)\,drdt\right|\\
\leq&\int^s_{s-\tau}(T-t)(\int^{T-t}_0(|\phi_t|^2+|\phi_r|^2+|h_t|^2)rdr)dt\\
\leq&(T-s+\tau)^2\tau.
\end{split}
\eeq
Letting $s\rightarrow T^-$ implies
\beq\label{term1}
\left|\iint_{K^T(\tau)}\left(r^2\phi_t\phi_r+r^2h_t\phi_r\right)\,drdt\right| \lesssim\tau^3.
\eeq
For second term, Green's theorem implies
\beq\notag
\begin{split}
&-\iint_{K^T(s,\tau)}\left((r^2\phi_t\phi_r)_t-\left(\frac{1}{2}r^2\phi_t^2+\frac{1}{2}r^2\phi_r^2-\frac{1}{2}\sin^2\phi\right)_r\right)\,drdt\\
=&\int^{T-s+\tau}_0r^2\phi_t\phi_r(r,s-\tau)\,dr+\int_{M^T(s,\tau)}r^2\phi_t\phi_r\,dr +\int_{T-s}^0r^2\phi_t\phi_r(r,s)\,dr\\
&+\int_{M^T(s,\tau)}\left(\frac{1}{2}r^2\phi_t^2+\frac{1}{2}r^2\phi_r^2-\frac{1}{2}\sin^2\phi\right)\,dt\\
=&\int^{T-s+\tau}_0r^2\phi_t\phi_r(r,s-\tau)\,dr+\int^0_{T-s}r^2\phi_t\phi_r(r,s)\,dr\\
&+\frac{1}{2}\int^s_{s-\tau}(\phi_t-\phi_r)^2(T-t,t)(T-t)^2\,dt-\frac{1}{2}\int^s_{s-\tau}\sin^2\phi(T-t,t)\,dt.
\end{split}
\eeq
Taking $s\rightarrow T^-$, we have
\beq\label{term2}
\begin{split}
&\left|\iint_{K^T(\tau)}\left((r^2\phi_t\phi_r)_t-\left(\frac{1}{2}r^2\phi_t^2+\frac{1}{2}r^2\phi_r^2-\frac{1}{2}\sin^2\phi\right)_r\right)\,drdt\right|\\
\lesssim& \int^\tau_0 e(r,T-\tau)r^2dr+\tau\text{Flux}(T,T-\tau)\\
=&\left(\int^{\delta\tau}_0+\int^\tau_{\delta\tau}\right)e(r,T-\tau)r^2dr+c\tau\text{Flux}(T,T-\tau)\\
\lesssim&\delta\tau E_0+\tau\int^\tau_{\delta\tau}e(r,T-\tau)rdr+\tau\text{Flux}(T,T-\tau)\\
\lesssim&  \epsilon\tau,
\end{split}
\eeq
 if we pick  $\delta=\frac{\epsilon}{3E_0}$ for any $\epsilon>0$ and use estimates \eqref{ext0} and \eqref{flux0}. Letting $s\rightarrow T^-$ in \eqref{pf2} and using \eqref{term1} and \eqref{term2}, we conclude \eqref{phit0}, which completes the proof of lemma.
\endpf

\section{Partial regularity and proof of Theorem \ref{thma}}
\setcounter{equation}{0}

In this section we will prove Theorem \ref{thma} on partial regularity. The proof will be divided into several lemmas. 

\begin{lemma}\label{englemma9}
Under the assumption of Theorem \ref{thma}, it holds
\beq
\sup_{(r,t)\in K^{T_0}_{\delta}}|\phi(r,t)|<c\epsilon_0,
\eeq
where
$$K^{T_0}_\delta=\{T_0-\delta\leq t\leq T_0,\ 0\leq r\leq T_0-t\}.$$
\end{lemma}

\pf
The proof is similar to Lemma \ref{englemma6}.
We may also need the following fact when $0<a<\pi$
$$
\int^a_0|\sin\psi|\,d\psi=\int^a_0\sin{\psi}\,d\psi\ =1-\cos{a}.
$$
This function is increasing w.r.t. $a$ and take maximum at $a=\pi$.
\endpf

\begin{remark} If we take $\epsilon_0$ small, then $\parallel\phi\parallel_{L^\infty}<\frac{\pi}{2}$ on $K^{T_0}_\delta$.
Therefore, there exists a positive constant $c$ such that
\beq\label{ineq2}
\frac{1}{c}\leq\frac{\sin^2\phi}{\phi^2}\leq c,
\eeq
and furthermore
\beq\label{ineq3}
\frac{1}{c}|\phi|\leq|\sin\phi||\cos\phi|\leq c|\phi|.
\eeq
\end{remark}

By the arguments of Lemma 3.1 in Struwe's paper \cite{struwe03}, we can show the following result. 
\begin{lemma}\label{englemma12}
Under the assumption of Theorem \ref{thma}, the solution $(\phi, h)$ to the system \eqref{welss} is H$\ddot{\mbox{o}}$lder continuous for any $r>0$ and $T>0$.
\end{lemma}

\begin{lemma}\label{englemma11}
Under the assumption of Theorem \ref{thma}, it holds
\beq\label{phicon}
\phi(r,t)\rightarrow 0\text{ as } (r,t)\rightarrow (0^+,T_0^-).
\eeq
which implies $\phi(r,t)$ is continuous at $(0,T_0)$.
\end{lemma}

\pf Denote $\chi(t)=\phi(T_0-t,t)$ for any $t>0$.
If we choose $\delta>0$ small enough in the assumption of Theorem \ref{thma}, then by Lemma \ref{englemma5},
$$
\text{Flux}\,(T_0,T_0-\delta)\rightarrow 0 \text{ as }\delta \rightarrow 0^+.
$$
By \eqref{ineq2} and the definition of Flux$\,(T_0,T_0-\delta)$, we obtain
$$
\int^{T_0}_{T_0-\delta}\frac{\chi^2(t)}{(T_0-t)^2}(T-t)\,dt\rightarrow 0 \text{ as } \delta\rightarrow 0^+,
$$
which implies that there exists a subsequence
$\delta_K\rightarrow 0^+$ such that $\delta_K<0$ and $\chi(T_0-\delta_K)\rightarrow 0$.
For any fixed $t\in (T_0-\delta,T_0)$,
\beq\notag
\begin{split}
&\left|\chi^2(t)-\chi^2(T_0-\delta_K)\right|\\
\lesssim& \int^t_{T_0-\delta_K}|\chi'||\chi|(s)\,ds\\
\lesssim&\left(\int^{T_0}_t|\chi'(s)|^2(T-s)\,ds\right)^\frac{1}{2}\left(\int^{T_0}_t\frac{|\chi(s)|^2}{(T-s)^2}(T-s)\,ds\right)^\frac{1}{2}\\
\lesssim&\ \text{Flux}\,(T_0,t),
\end{split}
\eeq
which implies
\beq\label{chi0}
\lim\limits_{t\rightarrow T_0^-}\chi(t)=\lim\limits_{t\rightarrow T_0^-}\phi(T_0-t,t)=0
\eeq

\medskip
Multiplying the second equation of \eqref{welss} by $tr\phi_t$, we obtain
\beq\label{3.7}
 \left(\frac{1}{2}tr\phi_t^2+\frac{1}{2}\frac{t\sin^2\phi}{r}+tr\frac{1}{2}\phi_r^2\right)_t-(tr\phi_r\phi_t)_r=
\frac{1}{2}r\phi_t^2+\frac{\sin^2\phi}{2r}+\frac{r}{2}\phi_r^2-tr\phi_t^2-trh_t\phi_t. 
\eeq
Integrating \eqref{3.7} over $K^T(s,s-\delta)$ for $0<s<T_0$, it holds
\beq\label{3.8}
\begin{split}
&\iint_{K^T(s,s-\delta)}[(tre(r,t))_t-(tr\phi_r\phi_t)_r]\,drdt\\
=&\iint_{K^T(s,s-\delta)}e(r,t)\,rdrdt-\iint_{K^T(s,s-\delta)}tr\phi_t^2+trh_t\phi_t\,drdt.
\end{split}
\eeq
By Green's theorem, the left side of \eqref{3.8} can be computed as follows
\beq\label{3.9}
\begin{split}
&\iint_{K^T(s,s-\delta)}[(tre(r,t))_t-(tr\phi_r\phi_t)_r]drdt
\\=&-\int^s_{s-\delta}t(T_0-t)\phi_r\phi_t(T_0-t,t)dt-\int^{T_0-(s-\delta)}_0(s-\delta)re(r,s-\delta)dr
\\&+\int^{T_0-s}_0sre(r,s)dr+\int^s_{s-\delta}t(T_0-t)e(T_0-t,t)dt.
\end{split}
\eeq
Letting $s \rightarrow T_0^-$ in \eqref{3.8} and \eqref{3.9}, we obtain
\beq\notag
\begin{split}
\int^\delta_0(T_0-\delta)e(r,T_0-\delta)\,rdr=&\int_{T_0-\delta}^{T_0}t(T_0-t)(e(T_0-t,t)-\phi_r\phi_t(T_0-t,t))\,drdt
\\&-\iint_{K^T(\delta)}e(r,t)\,rdrdt+\iint_{K^T(\delta)}tr\phi_t^2+trh_t\phi_t\,drdt,
\end{split}
\eeq
which implies
\beq\notag
\begin{split}
&\int^\delta_0e(r,T_0-\delta)\,rdr\\
\lesssim&\frac{T_0}{T_0-\delta}\,\text{Flux}\,(T_0,T_0-\delta)+\frac{1}{T_0-\delta}\iint_{K^T(\delta)}e(r,t)\,rdrdt\\
&+\frac{T_0}{T_0-\delta}\iint_{K^T(\delta)}(\phi_t^2+h_t^2)\,rdrdt.
\end{split}
\eeq
Combining this estimates with the energy estimates and Lemma \ref{englemma5}, we conclude 
\beq\label{simle}
\lim\limits_{t\rightarrow T_0^-}E(T_0-t,t)
=\lim\limits_{t\rightarrow T_0^-}\int_0^{T_0-t}e(r,t)\,rdr
= \lim\limits_{\delta\rightarrow 0^+}\int^\delta_0e(r,T_0-\delta)\,rdr
=0.
\eeq

\medskip
Now, we are ready to prove \eqref{phicon}. Direct computation implies
\beq\notag
\begin{split}
\phi^2(r,t)&=\phi^2(T_0-t,t)+\int^r_{T_0-t}2\phi\phi_r(r',t)\,dr'
\\&\lesssim \phi^2(T_0-t,t)+\left(\int^{T_0-t}_0\frac{\phi^2}{r^2}\,rdr\right)^\frac{1}{2}\left(\int_0^{T_0-t}\phi_r^2\,rdr\right)^\frac{1}{2}\\
&\lesssim\phi^2(T_0-t,t)+E(T_0-t,t). 
\end{split}
\eeq
Combining it with \eqref{chi0} and \eqref{simle}  completes the proof of lemma.
\endpf

\begin{lemma}\label{englemma13}
We have the following uniform estimate for the solution $h(r,t)$ to the system \eqref{welss} 
\beq\label{ineq4}
\parallel h_r\parallel_{L^\infty_{x,t}}+\parallel h_t\parallel_{L^\infty_{x,t}}\leq C.
\eeq
\end{lemma}

\begin{remark} Notice that here the smallness condition \eqref{ineq1} in Theorem \ref{thma} is not required. 
\end{remark}

\pf
We prove the conclusion by contradiction argument. Assume that $ h_r(r,t)\rightarrow \infty$ as $(r,t)\rightarrow (0^+,T_0^-)$. For any 
$ T_i\rightarrow T_0^-$, we 
denote
$$
\frac{1}{\lambda_i}=\sup_{0<r<\infty}h_r(r,T_i)\rightarrow \infty.
$$
Since $h_r$ has no singularity at $T_i$ and $h_r(\infty, T_i)=v(\infty, T_i)<\infty$, there exists $r_i\rightarrow 0^+$ such that 
$$
\frac{1}{\lambda_i}=\sup_{0<r<\infty}h_r(r,T_i)=h_r(r_i,T_i)
$$
For any $-\frac{r_i}{\lambda_i}<r<\infty$, $-\frac{T_i}{\lambda_i}\leq t\leq 0$, denote 
\beq\label{scahph}
h_{\lambda_i}(r,t)=h(r_i+\lambda_ir,T_i+\lambda_i^2t),\quad
\phi_{\lambda_i}(r,t)=\phi(r_i+\lambda_ir,T_i+\lambda_i^2t).
\eeq
Denote $T=T_i+\lambda_i^2t$, $R=r_i+\lambda_ir$ and direct computation implies that $h_{\lambda_i}$ satisfies the following equation
\beq\label{eqhli}
\begin{split}
(h_{\lambda_i})_t-(h_{\lambda_i})_{rr}
=\frac{(h_{\lambda_i})_r}{\frac{r_i}{\lambda_i}+r}-\frac{h_{\lambda_i}}{(\frac{r_i}{\lambda_i}+r)^2}+(\phi_{\lambda_i})_t
\end{split}
\eeq
By the choice of $\lambda_i$, it holds 
\begin{equation}\label{syseq1}
(h_{\lambda_i})_r(0,0)=1,\quad
\|(h_{\lambda_i})_r\|_{C^0}\leq 1.
\end{equation}
{
By Lemma \ref{englemma2} and invariance of the energy of $h$ under the scaling \eqref{scahph}, there exists a $\ h_\infty(r,t)$ such that  $h_{\lambda_i}\rightarrow h_\infty$ $W^*$ in $L^\infty_{\textrm{loc}}(\R,\dot{H})$ and $(h_{\lambda_i})_t\rightarrow (h_\infty)_t$ in $L^2_{\textrm{loc}}([0,\infty)\times(0,\infty))$ (up to subsequences, we still denote it as $h_{\lambda_i}$ for simplicity).
By Arzela-Ascoli theorem, there exists a function $\ h_\infty(r,t)$ such that $h_{\lambda_i}\rightarrow h_\infty \text{ in } C^0_{\textrm{loc}}[0,\infty)$ and 
\beq\label{hinft}
(h_\infty)_r(0,0)=1. 
\eeq
}
For any $N>0$ and $M>0$, direct computation implies
\beq\label{pfl2}
\begin{split}
\int_{-N}^0\int^M_{-M}|(h_{\lambda_i}(r,t))_t|^2\,rdrdt
=\int^{T_i}_{t_I-\lambda_i^2N}\int^{r_i+\lambda_iM}_{r_i-\lambda_iM}|h_T(R,T)|^2(R-r_i)\,dRdT\rightarrow 0\text{ as } i\rightarrow \infty,
\end{split}
\eeq
since $\lambda_i\rightarrow 0$, $r_i\rightarrow 0^+$, $T_i\rightarrow T_0^-$, which implies that $ h_\infty$ is independent of $t$.
Similarly, we have
\beq\label{pfl3}
\begin{split}
\int_{-N}^0\int^M_{-M}|(\phi_{\lambda_i}(r,t))_t|^2\,rdrdt
=\int^{T_i}_{t_I-\lambda_i^2N}\int^{r_i+\lambda_iM}_{r_i-\lambda_iM}|\phi_T(R,T)|^2(R-r_i)\,dRdT\rightarrow 0,
\end{split}
\eeq
and
\beq\label{pfl4}
\begin{split}
\int_{-M}^M|(h_{\lambda_i})_r|^2rdr&=\int_{r_i-\lambda_iM}^{r_i+\lambda_iM}|h_R(R,T)|^2RdR
<\infty.
\end{split}
\eeq

We will investigate the limit of $\frac{r_i}{\lambda_i}$ into three different cases. 

\medskip
\noindent{\em Case 1.} Assume that
\beq\label{cas1}
\lim_{i\rightarrow \infty}\frac{r_i}{\lambda_i}=\infty
\eeq
which implies the domain $h_{\infty}$ will be 
$(r,t)\in(-\infty,\infty)\times(-\infty,0)$.
By \eqref{syseq1}, we obtain
\beq\label{pfl1}
\frac{(h_{\lambda_i})_r}{\frac{r_i}{\lambda_i}+r}-\frac{h_{\lambda_i}}{(\frac{r_i}{\lambda_i}+r)^2}\rightarrow 0.
\eeq
Combining \eqref{pfl1} with \eqref{pfl2} and \eqref{pfl3}, the equation of $h_\infty$ would become
$$
(h_\infty)_{rr}=0 
$$
which can be solved directly by $h_\infty=C_1r+C_2$ for any constants $C_1$ and $C_2$. However, the estimate \eqref{pfl4} implies
$$
\int^\infty_{-\infty}|(h_\infty)_r|^2rdr<\infty
$$
which implies $C_1=0$ and then $(h_\infty)_r\equiv 0$. 
This contradicts to \eqref{hinft}.

\medskip
\noindent{\em Case 2.} Assume that
\beq\label{cas2}
\lim_{i\rightarrow \infty}\frac{r_i}{\lambda_i}=a>0,
\eeq
which implies the domain $h_{\infty}$ will be 
$(r,t)\in(-a,\infty)\times(-\infty,0)$. 
The equation of $h_\infty$ becomes
\beq\notag
-(h_\infty)_{rr}-\frac{(h_\infty)_r}{r+a}+\frac{h_\infty}{(r+a)^2}=0,\quad h_\infty(-a)=0.
\eeq
The general solution is $ h_\infty=C(r+a)$. By the fact $\int^\infty_{-a}|(h_\infty)_r|^2rdr<\infty$, we obtain $C\equiv 0$, which contradicts to \eqref{hinft}.


\medskip
\noindent{\em Case 3.} Assume that
\beq\label{cas3}
\lim_{i\rightarrow\infty}\frac{r_i}{\lambda_i}=0,
\eeq
which implies the domain $h_{\infty}$ will be 
$(r,t)\in(0,\infty)\times(-\infty,0)$. 
The equation of $h_\infty$ becomes 
\beq\notag
-(h_\infty)_{rr}-\frac{(h_\infty)_r}{r}+\frac{h_\infty}{r^2}=0,\quad
h_\infty(0)=0.
\eeq
The general solution is $h_\infty(r)=Cr$ for any constant $C$. The rest of proof is similar to last two cases, from which, we conclude that $\| h_r\|_{L^\infty_{r,t}}<\infty$.
%
%

\medskip
To estimate $\| h_t\|_{L^\infty_{r,t}}$, we may also assume by contradiction that $h_t(r,t)\rightarrow\infty$ as $(r,t)\rightarrow (0^+,T_0^-)$
Let $1/\lambda_i=\sup\limits_{0<r<\infty}\sqrt{h_t(r,T_i)}$
$\Rightarrow$ in the blow up sequence \eqref{scahph}. We obtain
$$
|(h_{\lambda_i}(r,t))_t|\leq 1,\quad (h_{\lambda_i})_t(0,0)=1.
$$
which contradicts to the fact that $h_{\infty}$ is independent of $t$. 
This completes the proof of lemma.
\endpf

\bigskip
\noindent{\bf Proof of Theorem \ref{thma}}. 
The partial regularity is now from the combination of Lemma \ref{englemma12}, Lemma \ref{englemma11} and Lemma \ref{englemma13}. 
\endpf

\section{Blowup and proof of Theorem \ref{thmb}}
\setcounter{equation}{0}

By Theorem \ref{thma}, we know if there is a blowup at $T_0$ iff 
\beq\label{4.1}
\inf_{0\leq t<T_0}E(T_0-t,t)\ge\epsilon_0>0.
\eeq
For a number $\epsilon_1>0$ ( will be determined later) and any $0\leq t<T_0$, we may choose $R=R(t)$ depends on $t$ such that
\beq\label{4.2}
\epsilon_1\leq E(6R,t)\leq 2\epsilon_1.
\eeq
By the local energy inequality \eqref{eneflux1}, we have
\beq\label{4.3}
E(R,t+\tau)\leq E(6R-\tau,t+\tau)\leq E(6R,t)+C\tau\leq 3\epsilon_1,
\eeq
for any $0<\tau\leq\min\{5R,\frac{\epsilon_1}{2C},\frac{1}{2}(T-t)\}$.
By similar arguments, we have
\beq\notag
\epsilon_1\leq E(6R,t)\leq E(6R+\tau,t-\tau)+C\tau\leq E(6R+\tau,t-\tau)+\frac{\epsilon_1}{2},
\eeq
which implies
\beq\label{4.4}
 E(6R+\tau,t-\tau)\ge\frac{\epsilon_1}{2}.
\eeq
By \eqref{4.1} $ and $ \eqref{4.3}, we should assume $3\epsilon_1<\epsilon_0$. Then we obtain $6R<T_0-t$ for and $ t\in[0,T_0)$.

\begin{lemma}\label{lemma4.1} According to the choice of $R$, it holds
\beq\label{4.5}
\lim_{t\rightarrow T_0^-}\frac{R}{T_0-t}=0.
\eeq
\end{lemma}

\pf
Suppose the limit is not valid. Then there exists sequences $t_i\rightarrow T_0^-$ and $R_i=R(t_i)$ such that 
$$
6R_i\ge C_0(T_0-t_i)
$$
for some constants $0<C_0<1$.
By assumption \eqref{4.1} $ and $ \eqref{4.2}, we have
$$
0<\epsilon_0-2\epsilon_1\leq E(T_0-t_i,t_i)-E(6R_i,t_i)\leq E(T_0-t_i,t_i)-E(C_0(T_0-t_i),t_i)
$$
which contradicts to Lemma \ref{englemma7}.
\endpf

\medskip
By Vitali's theorem, we can find countable disjoint intervals $I_i=(t_i-R_i,t_i+R_i)$ for $i=1, 2, \cdots$ such that $t_i$, $R_i$ satisfy the assumption \eqref{4.2} and $[0,T_0)\subset\bigcup\limits_{i=1}^\infty I_i^*$ where 
$
I_i^*=(t_i-5R_i,t_i+5R_i).
$
By the fact $6R_i<T_0-t_i$, it is not hard to see that 
\beq\label{4.8}
t_i+5R_i<T_0-R_i<T_0.
\eeq
Therefore, without loss of generality, we may assume $t_i\uparrow T_0^-$.

\begin{lemma}\label{lemma4.3}
It holds
\beq\label{4.6}
\lim_{i\rightarrow\infty}\frac{1}{R_i}\int_{I_i}\int^{T_0-t}_0|\phi_t|^2\,rdrdt= 0
\eeq
\end{lemma}

\pf
Assume that \eqref{4.6} is not true. Then $\exists\ \delta>0$ and $ i_0>0$ such that
\beq\label{4.7}
\int_{I_i}\int_0^{T_0-t}|\phi_t|^2rdrdt\ge \delta R_i
\eeq
for any $i\ge i_0$.
We obtain from \eqref{4.8} that
$\bigcup_{i<i_0}I_i^*$ can't cover $[0,T_0).$
For any $\tau\in(0,T_0)$ such that 
$T_0-\tau\in(\max\limits_{i<i_0}(t_i+5R_i),T_0)$
denote
$$
i_1=\max\{j:[T_0-\tau,T_0)\subset\bigcup_{i\ge j}I_i^*\}\ge i_0.
$$
The assumption \eqref{4.8} implies $T_0-\tau<t_{i_1}+5R_{i_1}<T_0-R_{i_1}$, which is equivalent to $R_{i_1}<\tau$.
Hence for any $i\ge i_1$, we obtain
$$
6R_i\leq T_0-t_i\leq T_0-t_{i_1}<\tau+5R_{i_1}<6\tau,
$$
which implies $t_i>T_0-6\tau$ and furthermore
$$
t_i-R_i>T_0-6\tau-\tau=T_0-7\tau,
$$
for any $i\ge i_1$. Then we obtain
\beq\label{4.9}
\begin{split}
\bigcup\limits_{i\ge i_1}I_i\subset [T_0-7\tau,T_0).
\end{split}
\eeq
Therefore by the definition of $i_1$, \eqref{4.7} and \eqref{4.8}, it holds
$$\delta\tau\leq\delta \sum_{i\ge i_1}\text{diam}(I_i^*)
=10\delta\sum_{i\ge i_1}R_i
\leq 10\sum_{i\ge i_1}\int_{I_i}\int_0^{T_0-t}|\phi_t|^2\,rdrdt
\leq 10\int^{T_0}_{T_0-7\tau}\int^{T_0-t}_0|\phi_t|^2\,rdrdt$$
which contradicts to lemma \ref{englemma8}. This completes the proof of lemma.
\endpf

\bigskip
\noindent{\bf Proof of Theorem \ref{thmb}}. \quad 
Denote
$$
\phi_i(r,t)=\phi(R_ir,t_i+R_it),\quad h_i(r,t)=h(R_ir,t_i+R_it).
$$
For any $t$ such that  $0<t_i+R_it<T_0$, it holds as $i\rightarrow \infty$
\beq\notag
\begin{split}
t<\frac{T_0-t_i}{R_i}\rightarrow \infty,\quad t>-\frac{t_i}{R_i}\rightarrow-\infty.
\end{split}
\eeq
Direct computation implies that 
$\phi_i$ and $h_i$ satisfy the following equations
\begin{equation}\label{4.10}
\begin{cases}
R_i(h_i)_t=\frac{1}{r}(r(h_i)_r)_r-\frac{h_i}{r^2}+R_i(\phi_i)_t\\
(\phi_i)_{tt}+R_i(\phi_i)_t=\frac{1}{r}(r(\phi_i)_r)_r-\frac{\sin(2\phi_i)}{2r^2}-R_i(h_i)_t.
\end{cases}
\end{equation}
Denote 
$$
r_i=\frac{T_0-t_i}{R_i}-1.
$$ 
It is not hard to see that $r_i\rightarrow \infty$
as $i\rightarrow\infty$. By Lemma \ref{lemma4.3}, we obtain
\beq\label{4.11}
\begin{split}
&\int_{-1}^1\int^{r_i}_0|(\phi_i)_t|^2\,rdrdt\\
=&\frac{1}{R_i}\int_{t_i-R_i}^{t_i+R_i}\int^{R_ir_i}_0|\phi_S(R,S)|^2\,RdRdS\\
\leq& \frac{1}{R_i}\int_{I_i}\int^{T-S}_0|\phi_S(R,S)|^2\,RdRdS\rightarrow 0
\end{split}
\eeq
where $R=R_ir$, $S=t_i+R_it$ and we have use the definition of $r_i$ or $R_ir_i=T_0-(t_i+R_i)$ in last inequality.
Similar to \eqref{pfl2}, for any $M,N>0$, we have
\beq\label{4.12}
\int_{-N}^N\int_0^M|R_i(h_i)_t|^2\,rdrdt
=R_i\int_{t_i-NR_i}^{t_i+NR_i}\int_0^{R_iM}|h_S|^2\,RdRdS\rightarrow 0,
\eeq
\beq\label{4.122}
\int_{-N}^N\int_0^M|R_i(\phi_i)_t|^2\,rdrdt
=R_i\int_{t_i-NR_i}^{t_i+NR_i}\int_0^{R_iM}|\phi_S|^2\,RdRdS\rightarrow 0.
\eeq
By the energy estimates, it holds for any time $t$
\beq\label{4.13}
\int_0^M|(\phi_i)_r|^2\,rdr=\int_0^{R_iM}|\phi_R(R,S)|^2\,RdR\leq E_0.
\eeq
Similarly
\beq\label{4.14}
\int_0^M|(h_{i})_r|^2rdr\leq E_0.
\eeq
Hence there are 
$\phi_\infty(r,t)$, $h_\infty(r,t)$ defined on $(r,t)\in(0,\infty)\times(-\infty,\infty)$ such that 
\beq\notag
\begin{split}
\phi_i\rightarrow \phi_\infty,\quad  h_i\rightarrow h_\infty \text{ weekly in } H^1_{loc}.
\end{split}
\eeq
By \eqref{4.11}, we know that $\phi_\infty$ is independent of $t$.
By \eqref{4.13}, \eqref{4.14} and weekly semi-continuity, we obtain
\beq\label{pf4.1}
\|(\phi_\infty)_{r}\|^2_{L^2(0,\infty)}+\|(h_\infty)_r\|^2_{L^2(0,\infty)}\leq E_0.
\eeq
By \eqref{4.12} and \eqref{4.122}, it is not hard to see that $\phi_\infty$ and $h_\infty$ are weak solution to the following system
\begin{equation}\label{4.15}\
\begin{cases}
h_{rr}+\frac{h_r}{r}-\frac{h}{r^2}=0,\\
\phi_{rr}+\frac{\phi_r}{r}-\frac{\sin\phi\cos\phi}{r^2}=0.
\end{cases}
\end{equation}
It is standard to show that $\phi_\infty$, $h_\infty$ are regular.
Solving for $h_\infty$ with $h_\infty(0,t)=0$, we have 
$h_\infty(r,t)=C(t)r$ or $h_\infty(r,t)=C(t)/r$. By the estimate \eqref{pf4.1}, it is not hard to see $h_\infty\equiv 0$.


Now we claim that $\phi_\infty$ is nontrivial solution to the second equqtion of \eqref{4.15} by using the fact \eqref{4.4} and showing $\phi_i\rightarrow \phi_\infty$ strongly in $H^1_{loc}$.
In fact, for any $\xi\in C_0^\infty(0,\infty)\times(-1,1)$ and $|\xi|\leq 1$.
Multiplying $(\phi_i-\phi_\infty)\xi$ for both equation of $\phi_i$ and $\phi_\infty$ and integrating over $(0,\infty)\times (-1,1)$, we have
\beq\label{4.16}
\begin{split}
&\iint|(\phi_i-\phi_\infty)_r|^2\xi \,rdrdt \\
=&\iint-(\phi_i)_{tt}(\phi_i-\phi_\infty)\xi\, rdrdt
-\iint R_i\phi_{it}(\phi_i\phi_\infty)\xi\, rdrdt
-\iint(h_i)_t(\phi_i-\phi_\infty)\xi\, rdrdt
\\&-\iint\left(\frac{\sin(2\phi_i)}{2r^2}-\frac{\sin(2\phi_\infty)}{2r^2}\right)(\phi_i-\phi_\infty)\xi \,r drdt
-\iint(\phi_i-\phi_\infty)_r(\phi_i-\phi_\infty)\xi_r\,rdrdt
\end{split}
\eeq
For the first term of right side, integrating by parts and using \eqref{4.11} and strong convergence of $\phi_i\rightarrow \phi_\infty$ in $L^2$, we obtain as $i\rightarrow \infty$
$$
\iint-(\phi_i)_{tt}(\phi_i-\phi_\infty)\xi\, rdrdt=\iint|(\phi_i)_t|^2\xi+(\phi_i)_t(\phi_i-\phi_\infty)\xi_t\,rdrdt\rightarrow 0.
$$
Similarly we can prove the second, third and fifth terms on right side approach to zero as $i\rightarrow \infty$.
The the last term, by \eqref{4.3} and the invariant of $L^2_r$ of $(\phi_i)_r$ and $(\phi_\infty)_r$, we have
$$
\|(\phi_i)_r\|_{L^2}^2,\ \|(\phi_\infty)_r\|_{L^2}^2<3\epsilon_1<\epsilon_0.
$$
By the proof of lemma \ref{englemma9},
$$
|\phi_i|, \ |\phi_\infty|<\pi. 
$$
Hence the last term of right side can be estimated as follows
$$
-\iint\left(\frac{\sin{2\phi_i}}{2r^2}-\frac{\sin(2\phi_\infty)}{2r^2}\right)(\phi_i-\phi_\infty)\xi\, rdrdt\leq 0.
$$
Combining all these facts with \eqref{4.16} we have 
$$
\phi_i\rightarrow\phi_\infty\text{ strongly in } H^1_{loc}.
$$
By the assumption \eqref{4.4}, we conclude that 
$\phi_\infty$ is a nontrival solution to the second equation of \eqref{4.15} and 
$\phi_\infty=2\arctan{\frac{r}{C}}$
for any constant $C$, which completes the proof of Theorem \ref{thmb}.
\endpf

\section{Proof of Theorem \ref{extthm}}
\setcounter{equation}{0}

The next lemma is similar to the proof of Lemma 1 in \cite{CKM}.

\begin{lemma}\label{lemma5.1} Suppose that $(\phi(r,t),h(r,t))$ is a solution to system \eqref{welss} with finite energy. Then for any $t>0$ we have
\beq\notag
h(t,0)=h(0,0),\quad h(t,\infty)=h(0,\infty),\quad
\phi(t,0)=\phi(0,0),\quad \phi(t,\infty)=\phi(0,\infty).
\eeq

\end{lemma}

\begin{lemma}\label{lemma5.2} Let $\phi(r,t)$ be a solution to system \eqref{welss} with
\beq\notag
\int \left(|\phi_r|^2+|\phi_t|^2+\frac{\sin^2\phi}{r^2}\right)\,rdr<4
\eeq
 and $\phi(0,t)=\phi(\infty,t)=0$. The following estimate is valid for any $(r,t)\in (0,\infty)\times(0,\infty)$
 \beq\label{6.1}
 |\phi(t,r)|\leq \frac{\pi}{2}.
 \eeq
\end{lemma}

\pf
Recall that $H(\phi(r,t))=\int_0^\phi|\sin\tau|d\tau$ in Lemma \ref{englemma7}. Since $H(\phi(0,t)=H(\phi(\infty,t)=0$, by the similar arguments as those in the proof of Lemma \ref{englemma7}, we should have 
\beq\notag
|H(\phi(t,r)|=|H(\phi(t,r))-H(\phi(t,0)|\leq\frac{1}{2}E(r,t)
\eeq
\beq\notag
|H(\phi(t,r)|=|H(\phi(t,r))-H(\phi(t,\infty)|\leq\frac{1}{2}\int_r^{\infty}\left(|\phi_r|^2+\frac{\sin^2\phi}{r^2}\right)\,rdr
\eeq
for any  $r>0$ and $t>0$. This implies
\beq\notag
|H(\phi(t,r))|\leq\frac{1}{4}\int_0^{\infty}\left(|\phi_r|^2+|\phi_t|^2+\frac{\sin^2\phi}{r^2}\right)\,rdr=1.
\eeq
Since H is a strictly increase function, we can conclude that 
\beq\notag
|\phi(r,t)|\leq H^{-1}(1)=\frac{\pi}{2}.
\eeq
\endpf

\bigskip
\noindent{\bf Proof of Theorem \ref{extthm}}. Suppose that $\phi$ blows up at some $0<T<\infty$. By Theorem \ref{thmb}, we can construct blowup sequence $\phi_i(r,t)$ such that 
\beq\notag
\phi_i(r,t)\rightarrow \phi_\infty=2\arctan\left(\frac{r}{C}\right) 
\eeq
in $H^1_{loc}$ for a nonconstant time-independent harmonic map $\phi_\infty(r)$.  By Lemma \ref{lemma5.2}, $\phi_\infty(r)\leq \frac{\pi}{2}$. However, by solving the equation of $\phi_\infty$ we obtain
$$\phi_\infty=2\arctan\left(\frac{r}{C}\right)$$
for any constant $C$, which is a contradiction to \eqref{6.1}.
Therefore, the solution $(\phi,h)$ has no blowup for any time $0<T<\infty$.
\endpf


%
\section*{Acknowledgments}
The first and third authors were partially supported by National Science Foundation (DMS-2306258).

\section*{Data Availability Statement}
This paper does not require any data.

\section*{Conflict-of-Interest Statement}
No potential conflict of interest was reported by the authors.

\end{document}